\documentclass[11pt]{article}
\usepackage{graphicx,wrapfig}
\graphicspath{figures/}
\usepackage{flafter}
\usepackage{leftidx}
\usepackage{mathrsfs}
\usepackage{amssymb,amsmath}
\usepackage{bbding}
\usepackage{fancyhdr}%
\usepackage{mathrsfs}
\usepackage{color}
\usepackage{stmaryrd}
\usepackage{amsfonts}
\usepackage{latexsym}
\usepackage{psfrag}
\usepackage{graphicx,subfigure}
\usepackage{fancyhdr,graphicx}
\usepackage{multicol}
\usepackage{dsfont}
\usepackage{bbm}
\usepackage{booktabs}
\usepackage[center]{caption2}
\usepackage{cite}
\usepackage{epstopdf}
\usepackage{diagbox}%
\usepackage{booktabs}%
\usepackage [latin1]{inputenc}%
\usepackage{multirow}%
\usepackage{enumerate}%
\usepackage{enumitem}
\usepackage{algpseudocode,algorithm,algorithmicx}%
\usepackage{appendix}%

\allowdisplaybreaks

\date{}
\newtheorem{theorem}{Theorem}[section]
\newtheorem{lemma}{Lemma}[section]
\newtheorem{remark}{Remark}
\newtheorem{example}{Example}[section]

\numberwithin{equation}{section}%

\newcommand{\rmk}{\mathrm{k}}
\newcommand{\defeq}{:=}
\newcommand{\zd}{\,\mathrm{d}}
\newcommand{\diff}{\triangledown_{\tau}}
\newcommand{\myvec}[1]{\boldsymbol{#1}}
\newcommand{\abs}[1]{\left|#1\right|}
\newcommand{\absb}[1]{\big|#1\big|}

\newcommand{\bra}[1]{\left(#1\right)}
\newcommand{\brab}[1]{\big(#1\big)}
\newcommand{\braB}[1]{\Big(#1\Big)}
\newcommand{\brat}[1]{(#1)}
\newcommand{\kbra}[1]{\left[#1\right]}
\newcommand{\kbrab}[1]{\big[#1\big]}

\newcommand{\myinner}[1]{\left\langle#1\right\rangle}
\newcommand{\myinnerb}[1]{\big\langle#1\big\rangle}

\newcommand{\mynorm}[1]{\left\|#1\right\|}
\newcommand{\mynormb}[1]{\big\|#1\big\|}

\newcommand{\timenorm}[1]{\absb{\!\absb{\!\absb{#1}\!}\!}}

\newcommand{\lan}[1]{{\color{blue}#1}}

\begin{document}
\title{$L^2$ norm error estimates of BDF methods up to
fifth-order for the phase field crystal model}
\author{Hong-lin Liao\thanks{ORCID 0000-0003-0777-6832; School of Mathematics,
Nanjing University of Aeronautics and Astronautics,
Nanjing 211106, China; Key Laboratory of Mathematical Modelling
and High Performance Computing of Air Vehicles (NUAA), MIIT, Nanjing 211106, China.
Hong-lin Liao (liaohl@nuaa.edu.cn and liaohl@csrc.ac.cn)
is supported by a grant 12071216 from National Natural Science Foundation of China.}
\quad
Yuanyuan Kang\thanks{School of Mathematics, Nanjing University of Aeronautics and Astronautics,
211101, P. R. China. Email: kangyy0101@163.com.}
}
\date{\today}
\maketitle
\normalsize

\begin{abstract}
  The well-known backward difference formulas (BDF) of the third, the fourth and the fifth orders
  are investigated for time integration of the phase field crystal  model.
  By building up novel discrete gradient structures of the BDF-$\rmk$ ($\rmk=3,4,5$) formulas,
  we establish energy dissipation laws at the discrete levels
  and then obtain a priori solution estimates for the associated numerical schemes
  (however, we can not build any discrete energy dissipation law for the corresponding
  BDF-6 scheme because the BDF-6 formula itself does not have any discrete gradient structures).
  With the help of the discrete orthogonal convolution kernels
  and Young-type convolution inequalities,
  some concise $L^2$ norm error estimates \lan{(with respect to the starting data
  	 in the $L^2$ norm)} 
  are established via the discrete energy technique.
  To the best of our knowledge, this is the first time such type $L^2$ norm error estimates of
  non-A-stable BDF schemes are obtained for nonlinear parabolic equations.
  Numerical examples are presented to verify and support the theoretical analysis.
  \noindent{\bf AMS subject classiffications.}\;\; 35Q99, 65M06, 65M12
\end{abstract}

\section{Introduction}\setcounter{equation}{0}

This work builds on the recent discrete energy analysis \cite{LiaoTangZhou:2021bdf345} of
the backward difference formula (BDF) schemes for linear diffusion equations.
The phase field crystal (PFC) model  is
a class of sixth order nonlinear parabolic equation, which is thermodynamically consistent
\cite{ElderKatakowskiHaatajaGrant:2002,ElderGrant:2004} in that the free energy
of the thermodynamic model is dissipative.
Consider a free energy functional of Swift-Hohenberg type \cite{ElderKatakowskiHaatajaGrant:2002,ElderGrant:2004},
\begin{align}\label{cont:free energy}
E[\Phi] = \int_{\Omega}\bra{\frac{1}{4}\Phi^4
+\frac{1}{2}\Phi\kbra{-\epsilon + (1+\Delta)^2}\Phi}\zd\mathbf{x},
\end{align}
where $\mathbf{x}\in\Omega\subseteq\mathbb{R}^2$,  $\Phi$ represents the
atomistic density field and $\epsilon\in(0,1)$ is a parameter related to the temperature.
The PFC equation is given by the $H^{-1}$ gradient flow
associated with the free energy functional \lan{$E[\Phi]$},
\begin{align}\label{cont: Problem-PFC}
\partial_t \Phi=\Delta\mu\quad\text{with the chemical potential}\quad
\mu:=\tfrac{\delta E}{\delta \Phi}= \Phi^3-\epsilon\Phi+(1+\Delta)^2\Phi.
\end{align}

The PFC growth model
is an efficient approach to
simulate crystal dynamics at the atomic scale in space while on diffusive scales in time.
This model has been successfully applied to
a wide variety of simulations in materials science across different
time scales. Related numerical schemes for the PFC model can be found
in \cite{DongFengWangWiseZhang:2018,HuWiseWangLowengrub:2009,
JingWang:2020,LiKim:2017, LiaoJiZhang:2020pfc,YangHan:2017}.
We assume that $\Phi$ is periodic over the domain $\Omega$.
Let the $L^2$ inner product $\myinner{f,g}:=\int_{\Omega}fg\zd{\mathbf{x}}$
and the associated $L^2$ norm $\mynorm{f}:=\sqrt{\myinner{f,f}}$ 
for all $f,g\in L^{2}(\Omega)$. We always use 
the standard seminorms and norms in the Sobolev space $H^{m}(\Omega)$ for $m\ge0$.
By the integration by parts, one has the volume conservation, $\myinnerb{\Phi(t),1}=\myinnerb{\Phi(0),1}$,
and the following energy dissipation law,
\begin{align}\label{cont:energy dissipation}
\frac{\zd{E}}{\zd{t}}=\brab{\tfrac{\delta E}{\delta \Phi},\partial_t \Phi}
=\bra{\mu,\Delta\mu}
=-\mynorm{\nabla\mu}^2 \le 0.
\end{align}

Let the discrete time level $t_k=k\tau$ with the uniform time-step
$\tau\defeq T/N$. For any discrete time sequence $\{v^n\}_{n=0}^N$,
denote $\diff v^n:=v^n-v^{n-1}$ and $\partial_{\tau}v^n:=\diff v^n/\tau.$
Here and hereafter, let the summation $\sum_{k=i}^{j}\cdot=0$
if the lower index $i$ is greater than the upper index $j$.
For a fixed index $3\le \rmk\le 6$, we view the BDF-$\rmk$  formula
as a discrete convolution summation,
\begin{align}\label{def: BDFk-Formula}
D_{\rmk}v^n:=\frac1{\tau}\sum_{j=1}^nb_{n-j}^{(\rmk)}\diff v^j
\quad \text{for $n\ge \rmk$},
\end{align}
where the associated  BDF-$\rmk$  kernels $b_{j}^{(\rmk)}$ (vanish if $j\ge \rmk$),
see Table \ref{table: BDF-k kernels}, are generated by
\begin{align}\label{def: BDF-k generating function}
\sum_{\ell=1}^{\rmk}\frac{1}{\ell}(1-\zeta)^{\ell-1}
=\sum_{\ell=0}^{\rmk-1}b_{\ell}^{(\rmk)}\zeta^\ell\quad\text{for $3\le \rmk\le 6$}.
\end{align}

\begin{table}[htb!]
\begin{center}\label{table: BDF-k kernels}
\caption{The BDF-$\rmk$ kernels $b_{j}^{(\rmk)}$ generated by \eqref{def: BDF-k generating function}}
\vspace*{0.3pt}
\def\temptablewidth{0.7\textwidth}
{\rule{\temptablewidth}{0.5pt}}
\begin{tabular*}{\temptablewidth}{@{\extracolsep{\fill}}ccccccc}
   BDF-$\rmk$      &$b_{0}^{(\rmk)}$     &$b_{1}^{(\rmk)}$     &$b_{2}^{(\rmk)}$
   &$b_{3}^{(\rmk)}$  &$b_{4}^{(\rmk)}$ &$b_{5}^{(\rmk)}$\\  \midrule
  $\rmk=3$   &$\frac{11}{6}$  &$-\frac{7}{6}$   &$\frac{1}{3}$   &&&\\[3pt]
  $\rmk=4$   &$\frac{25}{12}$ &$-\frac{23}{12}$ &$\frac{13}{12}$ &$-\frac{1}{4}$ &&\\[3pt]
  $\rmk=5$   &$\frac{137}{60}$&$-\frac{163}{60}$&$\frac{137}{60}$&$-\frac{21}{20}$ &$\frac{1}{5}$&\\[3pt]
  $\rmk=6$   &$\frac{147}{60}$&$-\frac{213}{60}$&$\frac{237}{60}$&$-\frac{163}{60}$ &$\frac{62}{60}$&$-\frac{1}{6}$\\
\end{tabular*}
{\rule{\temptablewidth}{0.5pt}}
\end{center}
\end{table}	

This work is motivated in developing high-order BDF time-steppings for the long-time
simulations of coarsening dynamics.
Recently, the adaptive BDF2 time-stepping scheme was investigated
theoretically in \cite{LiaoZhang:2021linear, LiaoJiZhang:2020pfc} for the linear diffusion equation and
PFC model \eqref{cont: Problem-PFC}, respectively.
The discrete energy dissipation law and concise $L^2$ norm error estimate
were established under a practical step-ratio constraint. In general, \lan{lower order schemes} would be well suited for
the fast varying solutions especially in the early coarsening process \cite{ZhangMaQiao:2013,LiLiao:2022ode,LiaoSongTangZhou:2020};
while high-order stable methods should be more preferable for slowly varying solutions
during the long-time process approaching the steady state \cite{BouchritiPierreAlaa:2020,HairerNorsettWanner:2002,HaoHuangWang:2021,Pierre:2021}.
\lan{It is well known that the BDF-$\rmk$ methods are numerically dissipative (L-stable) and have 
	a faster dissipation rate for higher frequency disturbances. 
	They are easy to implement compared with some existing methods and 
	have their own advantages in the long-time simulation of gradient flow problems, 
	including phase field crystal model.}
\lan{In the literature, the stability and convergence of A-stable (G-stable) BDF1 and BDF2 schemes \cite{Dahlquist:1978,ElliottStuart:1993,HillSuli:2000,HairerNorsettWanner:2002,StuartHumphries:1998}
	have been well studied, while the stability and convergence of the non-A-stable
BDF-$\rmk$ ($\rmk=3,4,5,6$) schemes for nonlinear phase field models
have not been well studied in the literatures due to the lack of proper discrete energy techniques.}

This situation was improved recently due to the seminal work \cite{LubichMansourVenkataraman:2013}
by Lubich, Mansour and Venkataraman. They noticed that the Nevanlinna-Odeh multiplier technique \cite{NevanlinnaOdeh:1981} is a powerful discrete tool for the stability analysis of non-A-stable BDF-$\rmk$ methods.
This tool was applied and explored in the numerical analysis of fully implicit and implicit-explicit
BDF-$\rmk$ approaches for linear and nonlinear parabolic problems, see related works in \cite{Akrivis:2015,AkrivisKatsoprinakis:2016,AkrivisLubich:2015} and references therein.
As noticed, the multiplier technique
relies on the celebrated equivalence of A-stability and G-stability
for linear multi-step methods by Dahlquist \cite{Dahlquist:1978}.
Recently, the discrete time derivative $D_{\rmk}v^n$ was also used in \cite{AkrivisFeischlKovacsLubich:2021,KovacsLiLubich:2019,KovacsLiLubich:2021} as a multiplier 
to derive optimal $H^1$ norm error estimate for nonlinear parabolic equations.
Nonetheless, because the products of nonlinear term and these multipliers
 can hardly be written into a difference between two positive functionals (part of energy), it seems that these multiplier techniques
are inadequate to establish the energy dissipation law and
$L^2$ norm convergence of non-A-stable
BDF-$\rmk$ schemes for nonlinear gradient flow problems.

Practically, the preservation of \eqref{cont:energy dissipation} at each time level,
called discrete energy dissipation law, has been proven to the fundamental requirement
of numerical methods for the effective simulation of long-time
coarsening dynamics
\cite{DongFengWangWiseZhang:2018,GongZhao:2019,HuWiseWangLowengrub:2009,
JingWang:2020,LiaoSongTangZhou:2020,Pierre:2021,StuartHumphries:1998}.
We focus on the intrinsic energy stability properties of
the non-A-stable BDF-$\rmk$ formulas themselves, that is,
some positive constants $\sigma_{L\rmk}$ (the larger, the better), two nonnegative
quadratic functionals $\mathcal{G}_{\rmk}$ and $\mathcal{R}_{\rmk}$ are sought such that
the BDF-$\rmk$ kernels $b_{j}^{(\rmk)}$ defined in \eqref{def: BDF-k generating function} satisfy
the following \textit{discrete gradient structure} in the sense of
\cite{BouchritiPierreAlaa:2020,Pierre:2021} or \cite[Section 5.6]{StuartHumphries:1998},
\begin{align}\label{def: discrete gradient structure}
v_n\sum_{j=1}^nb_{n-j}^{(\rmk)}v_j=&\,
\mathcal{G}_{\rmk}[\vec{v}_n]-\mathcal{G}_{\rmk}[\vec{v}_{n-1}]+\frac{\sigma_{L\rmk}}{2}v_n^2
+\mathcal{R}_{\rmk}[\vec{v}_n]\quad\text{for $n\ge \rmk$,}
\end{align}
where $\vec{v}_n$ denotes the consecutive tuples
$\vec{v}_n:=(v_n,v_{n-1},\cdots,v_{0})$.
As shown in Theorem \ref{thm:energy decay law} that the discrete gradient structure
\eqref{def: discrete gradient structure} plays an important role for constructing the
discrete energy dissipation laws of the corresponding BDF schemes.
In this work, we achieve concise discrete gradient structures for the BDF-$\rmk$ formulas
with the constants $\sigma_{L3}\approx1.979$, $\sigma_{L4}\approx1.601$ and $\sigma_{L5}\approx0.3367$, respectively,
see Lemma \ref{lem: Quadratic decomposition BDF-k}.

To demonstrate the practical significance of discrete gradient structures,  we consider
the following BDF-$\rmk$ implicit scheme subject to the periodic boundary conditions
\begin{align}\label{scheme: PFC BDFk}
D_{\rmk}\phi^n=\Delta\mu^n\quad\text{with}\quad
\mu^n=(1+\Delta)^2\phi^n + \brat{\phi^n}^3-\epsilon\phi^n
\quad\text{for $\rmk\le n\le N$,}
\end{align}
where the starting solutions $\phi^1$, $\phi^2$, $\cdots$, $\phi^{\rmk-1}$
are assumed to be available  and accurate enough, such as,  by Runge-Kutta methods \cite{GongZhao:2019}.
In this paper, we only consider the time-discrete approximation.
\lan{Our theoretical results including the discrete energy dissipation laws and $L^2$ norm error analysis 
	can be easily extended to the fully discrete scheme by using finite difference,
finite element or pseudo-spectral approximation preserving the discrete Green's formulas.}

In the $L^2$ norm error analysis, our main discrete tool is the
discrete orthogonal convolution (DOC) kernels.
For the discrete BDF-$\rmk$ kernels $b_{j}^{(\rmk)}$ generated by \eqref{def: BDF-k generating function},
the corresponding DOC-$\rmk$ kernels $\theta_{j}^{(\rmk)}$ are defined by \cite{LiaoTangZhou:2021bdf345}
\begin{align}\label{def: DOC-Kernels}
\theta_{0}^{(\rmk)}:=\frac{1}{b_{0}^{(\rmk)}}
\quad \mathrm{and} \quad
\theta_{n-j}^{(\rmk)}:=-\frac{1}{b_{0}^{(\rmk)}}
\sum_{\ell=j+1}^n\theta_{n-\ell}^{(\rmk)}b_{\ell-j}^{(\rmk)}\quad \text{for $j=n-1,n-2,\cdots, \rmk+1,\rmk$.}
\end{align}
It is easy to find the following \emph{discrete orthogonal convolution identity}
\begin{align}\label{eq: orthogonal identity}
\sum_{\ell=j}^{n}\theta_{n-\ell}^{(\rmk)}b^{(\rmk)}_{\ell-j}\equiv\delta_{nj}
\quad\text{for any $\rmk\leq j\le n$,}
\end{align}
where \lan{$\delta_{nj}$ is the Kronecker delta symbol.} Thus, by exchanging the summation order, one gets
\begin{align*}
\sum_{j=\rmk}^{n}\theta_{n-j}^{(\rmk)}
\sum_{\ell=\rmk}^{j}b_{j-\ell}^{(\rmk)}\diff \phi^\ell
=&\,\sum_{\ell=\rmk}^{n}\diff \phi^\ell\sum_{j=\ell}^{n}\theta_{n-j}^{(\rmk)}b_{j-\ell}^{(\rmk)}
=\diff \phi^n\quad\text{for $\rmk\le n\le N$.}
\end{align*}
Multiplying the BDF-k formula $D_{\rmk}\phi^j$ with the associated DOC kernels $\theta_{n-j}^{(\rmk)}$ 
and summing $j$ from $j=\rmk$ to $n$, we get
\begin{align}\label{Dis: DOC action BDF formula Dk}
\sum_{j=\rmk}^{n}\theta_{n-j}^{(\rmk)}D_{\rmk} \phi^j
=&\,\frac1{\tau}\sum_{j=\rmk}^{n}\theta_{n-j}^{(\rmk)}
\sum_{\ell=1}^{\rmk-1}b_{j-\ell}^{(\rmk)}\diff \phi^\ell
+\frac1{\tau}\sum_{j=\rmk}^{n}\theta_{n-j}^{(\rmk)}
\sum_{\ell=\rmk}^{j}b_{j-\ell}^{(\rmk)}\diff \phi^\ell\nonumber\\
\triangleq&\, \frac1{\tau}\phi_{\mathrm{I}}^{(\rmk,n)}+\partial_{\tau} \phi^n\qquad\text{for $\rmk\le n\le N$,}
\end{align}
where $ \phi_{\mathrm{I}}^{(\rmk,n)}$ represents the starting effects
on the numerical solution at the time $t_n$,
\begin{align}\label{Dis: initial effect -BDF formula}
 \phi_{\mathrm{I}}^{(\rmk,n)}:=\sum_{\ell=1}^{\rmk-1}\diff \phi^\ell
\sum_{j=\rmk}^{n}\theta_{n-j}^{(\rmk)}b_{j-\ell}^{(\rmk)}\qquad\text{for $n\ge\rmk$.}
\end{align}
By using \eqref{Dis: DOC action BDF formula Dk}
and \eqref{Dis: initial effect -BDF formula}, we can reformulate the discrete scheme \eqref{scheme: PFC BDFk} into
\begin{align}\label{eq: DOC form PFC BDFk}
\partial_{\tau} \phi^{j}=-\phi_{\mathrm{I}}^{(\rmk,j)}/\tau
+\sum_{\ell=\rmk}^j\theta_{j-\ell}^{(\rmk)}\Delta\mu^{\ell}\quad
\text{for $j\ge\rmk$}.
\end{align}
In section 4,  the $L^2$ norm error estimates of the BDF-$\rmk$ schemes \eqref{scheme: PFC BDFk} are proved
via the above equivalent formulation \eqref{eq: DOC form PFC BDFk}.
The standard discrete energy technique will be used with the help of 
some novel discrete convolution inequalities in section 3.
Numerical examples are presented in the last section to support our theoretical analysis.

In summary, our contributions in this paper are two-fold:
\begin{enumerate}
  \item Novel discrete gradient structures of the BDF-$\rmk$ ($\rmk=3,4,5$) formulas
are derived such that we can build up certain discrete energy dissipation laws
and obtain the priori solution estimates in the energy norm
for the BDF-$\rmk$ time-stepping schemes \eqref{scheme: PFC BDFk}.
However, we can not build any discrete energy dissipation law for the corresponding
  BDF-6 scheme because the BDF-6 formula itself does not have any discrete gradient structures.
  It provides a counterexample for the conjecture by Stuart and Humphries \cite[Section 5.6]{StuartHumphries:1998}.
  \item By developing novel discrete convolution inequalities of Young-type,
  \lan{we prove the $L^2$ norm convergence of the  high-order BDF-$\rmk$
scheme \eqref{scheme: PFC BDFk}  with respect to the starting data in the $L^2$ norm. To the best of our knowledge,
this is the first time such optimal $L^2$ norm error estimates of BDF-$\rmk$ methods ($\rmk=3,4,5$)
 are proved for a nonlinear parabolic problem.}
\end{enumerate}

Throughout this paper, any subscripted $C$, such as $C_\phi$,
denotes a generic positive constant, not necessarily
 the same at different occurrences; while,
 any subscripted $c$, such as $c_\Omega,c_0,c_1$ and so on,
denotes a fixed constant. Always, the appeared constants are  dependent on the given data
and the solution but independent of the time steps.

\section{Energy dissipation law and solvability}
\setcounter{equation}{0}

Denote the space 
$\mathbb{V}:=\{v\,|\,v\; \text{is periodic for}\; \mathbf{x}\in\Omega\}.$
For any functions $v,w\in\mathbb{V}$, one  has the Green's formulas,
$\myinner{-\Delta v,w}=\myinner{\nabla v,\nabla w}$,
$\myinner{\Delta^2v,w}=\myinner{\Delta v,\Delta w}$,  and $\myinner{\Delta^3v,w}=-\myinner{\nabla\Delta v,
	\nabla \Delta w}$.
Let $\mynorm{v}_{L^\infty}:=\max_{\mathbf{x}\in\Omega}|v|$. We have 
\lan{the embedding inequality \cite{ChengWangWiseYue:2016,LiaoSunShi:2010sinum}}
\begin{align}\label{ieq: Linfty H2 embedding}
	\mynormb{v}_{L^\infty}\le c_\Omega\bra{\mynormb{v}
		+\mynormb{\Delta v}}\quad \text{for any $v\in\mathbb{V}$.}
\end{align}
For the underlying volume-conservative problem,
it is convenient to define a mean-zero space
$\mathbb{\mathring V}:=\big\{v\in L^2(\Omega)\,|\, \myinner{v,1}=0\big\}\subset\mathbb{V}.$
The $H^{-1}$ inner product $\myinner{v,w}_{-1}:=\myinnerb{\bra{-\Delta}^{-1}v,w}$
and the associated $H^{-1}$ norm $\mynorm{\cdot}_{-1}$ can be defined by $\mynorm{v}_{-1}:=\sqrt{\myinner{v,v}_{-1}}\,.$
For any functions $v\in\mathbb{\mathring V}$, we have the generalized H\"{o}lder inequality, $\mynorm{v}^2\le \mynorm{\nabla v}\mynorm{v}_{-1}$, and the following lemma.
\begin{lemma}\cite[Lemma 2.1]{LiaoJiZhang:2020pfc}\label{lem: L2-Embedding-Inequality}
	For any grid functions $v\in \mathbb{\mathring V}$, it holds that
	$$\mynormb{v}^2 \le \frac{1}{3}\mynormb{(1+\Delta)v}^2+\frac{3}{2}\mynormb{v}_{-1}^2.$$
\end{lemma}

\subsection{Unique solvability}

To focus on the numerical analysis of the BDF-$\rmk$ solutions, it is to assume that
\begin{enumerate}[itemindent=1em]
\item[$\mathbf{A1}$.]
 Certain starting scheme, \lan{such as Gauss collocation Runge-Kutta 
 	method \cite{GongZhaoWang:2020HEQRK}}, is chosen to compute the first $(\rmk-1)$-level solutions
 $\phi^{\ell}$ for $1\le \ell\le \rmk-1$ such that they preserve the volume,
 $\myinnerb{\phi^\ell,1}=\myinnerb{\phi^0,1}$ for $1\le \ell\le \rmk-1$.
  \end{enumerate}

Note that, the solution $\phi^n$ of BDF-$\rmk$ scheme \eqref{scheme: PFC BDFk}
preserves the volume, $\myinnerb{\phi^n,1}=\myinnerb{\phi^0,1}$,
for $n\ge1$. Actually, taking the inner product of \eqref{scheme: PFC BDFk} by 1 and applying the summation by parts, one has
 $\myinnerb{D_{\rmk}\phi^j,1}=\myinnerb{\Delta\mu^j,1}=0$ for $j\ge \rmk$.
 Multiplying both sides of this equality by the DOC-$\rmk$ kernels
 $\theta_{n-j}^{(\rmk)}$ and summing the index $j$ from $j=\rmk$ to $n$, we get
 $\sum_{j=\rmk}^n\theta_{n-j}^{(\rmk)}\myinnerb{D_{\rmk}\phi^j,1}=0$ for $n\ge\rmk$.
It follows from \eqref{Dis: DOC action BDF formula Dk} that $\myinnerb{\diff \phi^n,1}=0$
because the assumption $\mathbf{A1}$ implies
$\myinnerb{\phi_{\mathrm{I}}^{(\rmk,n)},1}=0$.
Simple induction yields the conservation law, $\myinnerb{\phi^n,1}=\myinnerb{\phi^{0},1}$
for $n\ge1$.

\begin{theorem}\label{thm: convexity solvability}
If  the time-step $\tau\le \frac2{3\epsilon}b_0^{(\rmk)}$,
the BDF-$\rmk$ scheme \eqref{scheme: PFC BDFk} is uniquely solvable.
\end{theorem}

\begin{proof}
For any fixed time-level indexes $n\ge\rmk$,
we consider the following energy functional $G$ on the space
$\mathbb{V}^{*}:=\big\{z\in L^2(\Omega)\,|\, \myinnerb{z,1}=\myinnerb{\phi^{n-1},1}\big\},$
\begin{align*}
G[z]:=&\,\frac1{2\tau}\myinnerb{b_0^{(\rmk)}\brat{z-\phi^{n-1}}+2L^{n-1},z-\phi^{n-1}}_{-1}
+\frac12\mynormb{(1+\Delta)z}^2+\frac14\myinnerb{z^3-2\epsilon z,z},
\end{align*}
where
$L^{n-1}:=\sum_{\ell=1}^{n-1}b_{n-\ell}^{(\rmk)}\diff \phi^{\ell}.$
Under the time-step constraint $\tau\le \frac2{3\epsilon}b_0^{(\rmk)}$, the functional $G$ is strictly convex.
Actually, for any $\lambda\in \mathbb{R}$
and any $\psi\in \mathbb{\mathring V}$, one has
\begin{align*}
\frac{\zd^2G}{\zd\lambda^2}[z+\lambda\psi]\Big|_{\lambda=0}
=&\,\frac1{\tau}b_0^{(\rmk)}\mynormb{\psi}_{-1}^2
+\mynormb{(1+\Delta)\psi}^2
+3\mynormb{z\psi}^2-\epsilon\mynormb{\psi}^2\\
\ge&\,\brab{\frac1{\tau}b_0^{(\rmk)}-\frac{3\epsilon}2}\mynormb{\psi}_{-1}^2
+\frac{2}3\mynormb{(1+\Delta)\psi}^2
+3\mynormb{z\psi}^2\ge0,
\end{align*}
where Lemma \ref{lem: L2-Embedding-Inequality} was applied with the setting $0<\epsilon<1$.
 Also, $G[z]$ is coercive on $\mathbb{V}^{*}$.
Thus the functional $G$ has a unique minimizer, denoted by $\phi^n$, if and only if it solves the equation
\begin{align*}
0=&\,\frac{\zd G}{\zd\lambda}[z+\lambda\psi]\Big|_{\lambda=0}
=\frac1{\tau}\myinnerb{b_0^{(\rmk)}\brat{z-\phi^{n-1}}+L^{n-1},\psi}_{-1}
+\myinnerb{(1+\Delta)^2z+z^3-\epsilon z,\psi}\\
=&\,
\frac1{\tau}\myinnerb{b_0^{(\rmk)}\brat{z-\phi^{n-1}}+L^{n-1}
-\tau\Delta\kbrab{(1+\Delta)^2z+z^3-\epsilon z},\psi}_{-1}.
\end{align*}
This equation holds for any $\psi\in \mathbb{\mathring V}$ if and only
if the unique minimizer $\phi^n\in\mathbb{V}^{*}$ solves
\begin{align*}
\frac1{\tau}\sum_{\ell=1}^{n}b_{n-\ell}^{(\rmk)}\diff \phi^{\ell}
-\Delta\kbrab{(1+\Delta)^2\phi^n+(\phi^n)^3-\epsilon \phi^n}=0,
\end{align*}
which is just the BDF-$\rmk$ scheme \eqref{scheme: PFC BDFk}.
It  completes the proof.
\end{proof}

\lan{As seen, the time-step size constraint	$\tau\le\frac{2}{3\epsilon}b_0^{(\rmk)}$ for solvability
	and energy stability (see Theorem \ref{thm:energy decay law} below)
is not practically restrictive since the parameter $\epsilon\in(0,1)$.}

\subsection{Energy dissipation law}

The positive definiteness of the BDF-$\rmk$ kernels has been established
in \cite[Lemma 2.4]{LiaoTangZhou:2021bdf345}
with the help of the Grenander-Szeg\"{o} theorem \cite[pp. 64--65]{GrenanderSzego:2001}.

\begin{lemma}\cite[Lemma 2.4]{LiaoTangZhou:2021bdf345}\label{lem: BDF-k minimum eigenvalue}
For $3\le \rmk\le 5$,
the discrete BDF-$\rmk$ kernels $b_{j}^{(\rmk)}$ defined in \eqref{def: BDF-k generating function}
are positive definite in the sense that
\[
2\sum_{\ell=\rmk}^n w_{\ell} \sum_{j=\rmk}^\ell b_{\ell-j}^{(\rmk)}w_j
\ge \mathfrak{m}_{1\rmk}\sum_{\ell=\rmk}^nw_\ell^2\quad\text{for $n\ge \rmk$},
\]
where $\mathfrak{m}_{13}=95/48$, $\mathfrak{m}_{14}=1.628$ and $\mathfrak{m}_{15}=0.3711$.
The minimum eigenvalue of the associated quadratic form with the BDF-$\rmk$ kernels
can be bounded from below by the constant $\mathfrak{m}_{1\rmk}$.
\end{lemma}

This result may be adequate to show that the discrete solution of \eqref{scheme: PFC BDFk}
is bounded in an energy norm. However, it is inadequate to build
some discrete energy dissipation laws to simulate the continuous property \eqref{cont:energy dissipation}
at each time level.
To achieve this aim, some novel quadratic decompositions
(or, discrete gradient structures according to
\cite{ElliottStuart:1993,StuartHumphries:1998}) for
the BDF-$\rmk$ formulas \eqref{def: BDFk-Formula} are given in the following lemma.
Some roughly lower estimates $\sigma_{L\rmk}\le \mathfrak{m}_{1\rmk}$
are then obtained for the minimum eigenvalues of the quadratic forms with the BDF-$\rmk$ kernels $b_{j}^{(\rmk)}$.


We remark that the functionals $\mathcal{G}_{\rmk}$ and $\mathcal{R}_{\rmk}$ 
in Lemma \ref{lem: Quadratic decomposition BDF-k}
always involve the consecutive tuples $\vec{v}_n=(v_n,v_{n-1},\cdots,v_{0})$.  
For the simplicity of notations, we denote
\begin{align*}
	\mathcal{G}_{\rmk}[\vec{v}_n]\triangleq\mathcal{G}_{\rmk}[v_n,v_{n-1},\cdots,v_{0}],\quad
	\mathcal{G}_{\rmk}[\vec{v}_{n-1}]\triangleq\mathcal{G}_{\rmk}[v_{n-1},v_{n-2},\cdots,v_{0}].
\end{align*}

\begin{lemma}\label{lem: Quadratic decomposition BDF-k}
For the real sequence $\{v_k\,|\,k=0,1,2,\cdots, N\}$,
define the difference operators
$$\delta_1v_n:=\delta_1^1v_n=v_n-v_{n-1}\quad
\text{and}\quad \delta_1^{m+1}v_n:=\delta_1^{m}(\delta_1v_n)
=\delta_1^{m}v_n-\delta_1^{m}v_{n-1}\quad \text{for $m\ge1$.}$$
Then for the step index $\rmk=3,4$ and 5, there exists positive constant $\sigma_{L\rmk}$,
nonnegative quadratic functionals $\mathcal{G}_{\rmk}$ and $\mathcal{R}_{\rmk}$ such that
 the BDF-$\rmk$ kernels $b_{j}^{(\rmk)}$ defined in \eqref{def: BDF-k generating function} satisfy
\begin{align}\label{QuadDecompo: general quadratic decomposition}
v_n\sum_{j=1}^nb_{n-j}^{(\rmk)}v_j=&\,
\mathcal{G}_{\rmk}[\vec{v}_n]-\mathcal{G}_{\rmk}[\vec{v}_{n-1}]
+\frac{\sigma_{L\rmk}}{2}v_n^2+\mathcal{R}_{\rmk}[\vec{v}_n]\quad\text{for $n\ge \rmk$, }
\end{align}
where positive constants $\sigma_{L\rmk}$, the quadratic functionals
 $\mathcal{G}_{\rmk}$ and $\mathcal{R}_{\rmk}$ are given by
\begin{itemize}
  \item for $\rmk=3$, the constant $\sigma_{L3}:=\frac{95}{48}\approx1.979$,
  \begin{align*}
\mathcal{G}_{3}[\vec{v}_n]:=&\,\frac{37}{96}v_{n}^2-\frac{1}{8}v_{n-1}^2+\frac{7}{24}(\delta_1v_{n})^2
=\frac{1}{6}v_{n}^2+\frac1{6}(\tfrac7{4}v_n-v_{n-1})^2,\\
\mathcal{R}_{3}[\vec{v}_n]:=&\,\frac{1}{6}(\delta_1^2v_n+\tfrac14v_{n-1})^2;\hspace{8cm}
\end{align*}
   \item for $\rmk=4$, the constant $\sigma_{L4}:=\frac{4919}{3072}\approx1.601$,
     \begin{align*}
\mathcal{G}_{4}[\vec{v}_n]:=&\,\frac{3433}{6144}v_{n}^2-\frac{15}{64}v_{n-1}^2+\frac{1}{8}v_{n-2}^2
+\frac{47}{192}(\delta_1v_{n})^2-\frac{3}{16}(\delta_1v_{n-1})^2+\frac{3}{16}(\delta_1^2v_{n})^2\\
=&\,\frac{13627}{43008}v_{n}^2+\frac{7}{24}(\tfrac{65}{56}v_n-v_{n-1})^2
+\frac{1}{8}(\tfrac3{2}\delta_1v_{n}+v_{n-2})^2,\\
\mathcal{R}_{4}[\vec{v}_n]:=&\,\frac{1}{8}(\delta_1^3v_n+\tfrac3{2}\delta_1v_{n-1})^2
+\frac{1}{6}(\delta_1^2v_n+\tfrac{35}{32}v_{n-1})^2;
\end{align*}
  \item for $\rmk=5$, the constant $\sigma_{L5}:=\frac{646631}{1920000}\approx0.3367$,
       \begin{align*}
\mathcal{G}_{5}[\vec{v}_n]:=&\,\frac{4227769}{3840000}v_{n}^2-\frac{551}{1600}v_{n-1}^2+\frac{17}{40}v_{n-2}^2
-\frac{1}{10}v_{n-3}^2+\frac{1607}{4800}(\delta_1v_{n})^2\nonumber\\
&\,-\frac{39}{80}(\delta_1v_{n-1})^2+\frac{2}{5}(\delta_1v_{n-2})^2
+\frac{7}{80}(\delta_1^2v_{n})^2-\frac{2}{5}(\delta_1^2v_{n-1})^2
+\frac{1}{5}(\delta_1^3v_{n})^2\\
=&\,\frac{1198850903}{1678080000}v_{n}^2
+\frac{437}{900}(\tfrac{4931}{6992}v_{n}-v_{n-1})^2\nonumber\\
&\,+\frac{9}{40}(\tfrac{23}{18}\delta_1v_{n}+v_{n-2})^2
+\frac{1}{10}(2\delta_1v_{n}+2v_{n-2}-v_{n-3})^2,\\
\mathcal{R}_{5}[\vec{v}_n]:=&\,\frac{1}{10}(\delta_1^4v_n+2\delta_1^2v_{n-1})^2
+\frac{1}{8}(\delta_1^3v_n+\tfrac{23}{10}\delta_1v_{n-1})^2
+\frac{1}{6}(\delta_1^2v_n+\tfrac{1787}{800}v_{n-1})^2.
\end{align*}
\end{itemize}
Then, by summing \eqref{QuadDecompo: general quadratic decomposition} 
with the setting $v_j=0$ for $1\le j \le \rmk-1$, 
 the associated quadratic form of BDF-$\rmk$ kernels $b_{j}^{(\rmk)}$ can be bounded by
\begin{align*}
2\sum_{\ell=\rmk}^nv_\ell\sum_{j=\rmk}^\ell b_{\ell-j}^{(\rmk)}v_j
\ge\sigma_{L\rmk}\sum_{\ell=\rmk}^nv_\ell^2\quad\text{for $n\ge \rmk$.}
\end{align*}
\end{lemma}

The discrete gradient structures \eqref{QuadDecompo: general quadratic decomposition}
for the BDF-$\rmk$ formulas with $\rmk=3,4,5$  can be checked
by some symbolic computation software
 (see the appended MATHEMATICA program \verb"Appendix_BDF345decomposition.nb")
 or by rather lengthy but delicate calculations
(Appendix \ref{section: quadratic decomposition BDF-k} gives a detail proof
of Lemma \ref{lem: Quadratic decomposition BDF-k} for interested readers).
Note that the quadratic decomposition
for the case of $\rmk=3$ would be optimal in the sense that
the resulting minimum eigenvalue bound $\sigma_{L3}=95/48$ equals the lower bound $\mathfrak{m}_{13}$,
see Lemma \ref{lem: BDF-k minimum eigenvalue}.
The cases of $\rmk=4,5$ seem to be nearly optimal in the sense that
$\sigma_{L4}$ and $\sigma_{L5}$ is very close to $\mathfrak{m}_{14}$ and $\mathfrak{m}_{15}$, respectively. 

The delicate quadratic decompositions \eqref{QuadDecompo: general quadratic decomposition}
for the non-A-stable BDF-$\rmk$ methods significantly update the results
in \cite[Theorem 5.6.3]{StuartHumphries:1998} or \cite[Theorem 6.2]{ElliottStuart:1993}.
They are \lan{much sharper than} the recent results in \cite[Theorems 3.2 and 3.6]{BouchritiPierreAlaa:2020}
with the eigenvalue estimates $\sigma_{L4}= 41/72\approx0.5694$ and $\sigma_{L5}= 0.1$
for the BDF-4 and BDF-5 formulas, respectively. The present results give the explicit
expressions of the Lyapunov functionals $\mathcal{G}_{\rmk}$
 and the proof is quite different from the technique of undetermined coefficients in
 \cite{BouchritiPierreAlaa:2020,Pierre:2021}.

\begin{remark}\label{remark: BDF-6 scheme}
The BDF-6 formula might not be suited for simulating the gradient flow models,
because we can not find two nonnegative quadratic functionals
 $\mathcal{G}_{6}$ and $\mathcal{R}_{6}$ to ensure the discrete gradient structure
\eqref{QuadDecompo: general quadratic decomposition} for the BDF-6 formula.
Otherwise, the discrete BDF-6 kernels $b_{j}^{(6)}$ defined by \eqref{def: BDF-k generating function}
are at least positive semi-definite. However, it is not difficult to check that the associated
 quadratic form $\sum_{m=6}^nv_m\sum_{j=6}^mb_{m-j}^{(6)}v_j$ has negative eigenvalues for proper large $n$.
Moreover, the BDF-6 formula is A$(\alpha)$-stable with $\alpha=17.84^\circ$,
cf. \cite[Section V.2]{HairerNorsettWanner:2002}.
Thus it provides a counterexample for the conjecture by Stuart and Humphries in
\cite[Section 5.6]{StuartHumphries:1998}, in which they inferred that
 A$(\alpha)$-stability implies  ``gradient stability". 
\end{remark}


Let $E[\phi^n]$ be the discrete version of free energy functional \eqref{cont:free energy}, given by
\begin{align}\label{def: discrete free energy}
E[\phi^n]
:=\frac{1}{2}\mynormb{(1+\Delta)\phi^n}^2 +
\frac{1}{4}\mynormb{(\phi^n)^2-\epsilon}^2-\frac{1}{4}\mynormb{\epsilon}^2 \quad\text{for $n\ge 0$.}
\end{align}
We define a modified discrete energy for $n\ge \rmk$,
\begin{align}\label{def: modified discrete energy}
\mathcal{E}_{\rmk}[\vec{\phi}^n]
:=&\,E[\phi^n] + \frac{1}{\tau}\myinnerb{\mathcal{G}_{\rmk}\kbrab{\diff\vec{\phi}^n},1}_{-1}.
\end{align}
As seen, the modified discrete energy $\mathcal{E}_{\rmk}$ introduces a perturbed term
of $O(\tau)$ to the original energy $E[\phi^n]$ due to the application of BDF-$\rmk$ formula $D_{\rmk}$.


\begin{theorem}\label{thm:energy decay law}
Assume that \textbf{A1} holds and the time-step sizes are properly small such that
\begin{align}\label{Restriction-Time-Step}
\tau\le\frac{2}{3\epsilon}\min\big\{b_0^{(\rmk)},\sigma_{L\rmk}\big\},
\end{align}
where $\sigma_{L3}\approx1.979>b_{0}^{(3)}$,
$\sigma_{L4}\approx1.601<b_{0}^{(4)}$
and $\sigma_{L5}\approx0.3367<b_{0}^{(5)}$.
Then the BDF-$\rmk$ implicit scheme \eqref{scheme: PFC BDFk}
preserves the following energy dissipation law
\begin{align*}
\mathcal{E}_{\rmk}[\vec{\phi}^n] \le \mathcal{E}_{\rmk}[\vec{\phi}^{n-1}]
\quad\text{for $n\ge \rmk$.}
\end{align*}
\end{theorem}
\begin{proof}
The first bound of \eqref{Restriction-Time-Step} ensures the solvability
in Theorem \ref{thm: convexity solvability}.
We shall establish the energy law under the second condition of \eqref{Restriction-Time-Step}.
The volume conservation law implies $\diff \phi^n\in\mathbb{\mathring V}$ for $n\ge1$.
Then we make the inner product of \eqref{scheme: PFC BDFk} by $(-\Delta)^{-1}\diff \phi^{n}$ and obtain
\begin{align}\label{Energy-Law-Inner}
\myinnerb{D_\rmk\phi^n,(-\Delta)^{-1}\diff \phi^{n}}
+\myinnerb{(1+\Delta)^2\phi^n,\diff \phi^{n}}
+\myinnerb{\brat{\phi^n}^3-\epsilon \phi^n,\diff \phi^{n}}=0.
\end{align}
With the help of the summation by parts and $2a(a-b)=a^2-b^2+(a-b)^2$,
the second term at the left hand side of \eqref{Energy-Law-Inner} gives
\begin{align*}
\myinnerb{(1+\Delta)^2\phi^n,&\,\diff \phi^n}
=\frac{1}{2}\mynormb{(1+\Delta)\phi^n}^2
-\frac{1}{2}\mynormb{(1+\Delta)\phi^{n-1}}^2
+\frac{1}{2}\mynormb{(1+\Delta)\diff \phi^{n}}^2.
\end{align*}
By using Lemma \ref{lem: Quadratic decomposition BDF-k} with $v_n:=\diff \phi^n$,
the first term in \eqref{Energy-Law-Inner} can be bounded by
\begin{align*}
\myinnerb{D_{\rmk}\phi^n,(-\Delta)^{-1}\diff \phi^{n}}
\ge&\, \frac{1}{\tau}\myinnerb{\mathcal{G}_{\rmk}\kbrab{\diff\vec{\phi}^n},1}_{-1}
-\frac{1}{\tau}\myinnerb{\mathcal{G}_{\rmk}\kbrab{\diff\vec{\phi}^{n-1}},1}_{-1}
+\frac{\sigma_{L\rmk}}{2\tau}\mynormb{\diff \phi^n}_{-1}^2.
\end{align*}
It is easy to check the following identity
\begin{align*}
4\brab{a^3-\epsilon a}\bra{a-b}= \bra{a^2-\epsilon}^2 - \bra{b^2-\epsilon}^2-2\bra{\epsilon-a^2}\bra{a-b}^2 + \bra{a^2-b^2}^2.
\end{align*}
One can bound the third term in \eqref{Energy-Law-Inner} by
\begin{align*}
\myinnerb{(\phi^n)^3-\epsilon\phi^n,\diff \phi^n}
\ge&\, \frac{1}{4}\mynormb{(\phi^n)^2-\epsilon}^2 - \frac{1}{4}\mynormb{(\phi^{n-1})^2-\epsilon}^2
-\frac{\epsilon}{2}\mynormb{\diff \phi^n}^2.
\end{align*}
By collecting the above estimates, it follows from \eqref{Energy-Law-Inner}
and the definition \eqref{def: modified discrete energy} that
\begin{align}\label{Energy-Inequality}
\frac{1}{2}\mynormb{(1+\Delta)\diff \phi^{n}}^2
+\frac{\sigma_{L\rmk}}{2\tau}\mynormb{\diff \phi^n}_{-1}^2
-\frac{\epsilon}{2}\mynormb{\diff \phi^n}^2
+\mathcal{E}_{\rmk}[\vec{\phi}^n]\le \mathcal{E}_{\rmk}[\vec{\phi}^{n-1}]
\end{align}
for $n\ge \rmk$.
Applying Lemma \ref{lem: L2-Embedding-Inequality}, one has
\begin{align*}
\frac{\epsilon}{2}\mynormb{\diff \phi^n}^2 \le \frac{1}{6}\mynormb{(1+\Delta)\diff \phi^n}^2
+\frac{3\epsilon}{4}\mynormb{\diff \phi^n}_{-1}^2,
\end{align*}
where $0<\epsilon<1$ has been used.
Thus we can obtain that
\begin{align*}
\frac{1}{3}\mynormb{(1+\Delta)\diff \phi^{n}}^2
+\brab{\frac{\sigma_{L\rmk}}{2\tau}-\frac{3\epsilon}{4}}\mynormb{\diff \phi^n}_{-1}^2+\mathcal{E}_{\rmk}[\vec{\phi}^n]\le \mathcal{E}_{\rmk}[\vec{\phi}^{n-1}]\quad\text{for $n\ge \rmk$.}
\end{align*}
Under the second condition of \eqref{Restriction-Time-Step}
or $\tau\le \frac{2}{3\epsilon}\sigma_{L\rmk}$, it yields the claimed result.
\end{proof}

The two time-step constraints in \eqref{Restriction-Time-Step} ensure the unique solvability
and the energy stability are consistent since they have the same order of magnitude.
But the constraint \eqref{Restriction-Time-Step} always requires smaller step-sizes
for the higher order methods. It is expected that some stabilized techniques \cite{XuTang:2006}
 would remove the time-step restriction without sacrificing the time accuracy. However,
this issue is out of our current scope
and will be reported in further studies.


To simplify the subsequent analysis, we impose a further assumption:
\begin{enumerate}[itemindent=1em]
\item[$\mathbf{A2}$.]
  Under the assumption $\mathbf{A1}$ and the time-step constraint \eqref{Restriction-Time-Step},
  assume that there exists a constant $c_0$ such that $\mathcal{E}_{\rmk}[\vec{\phi}^{\rmk-1}]\le c_0$,
  where $c_0$ may depend on the problem and the starting values,
  but is always independent of the time-step size $\tau$.
  \end{enumerate}

\begin{lemma}\label{lem:Bound-Solution}
If $\mathbf{A2}$ holds, the solution of BDF-$\rmk$ scheme
\eqref{scheme: PFC BDFk} is stable in the $L^{\infty}$ norm,
$$\mynormb{\phi^n}_{\infty}\le
c_1:=c_{\Omega}\sqrt{8c_0+2\brab{2+\epsilon}^2\abs{\Omega}}\quad\text{for $n\ge \rmk$},$$
where $c_1$ may depend on the problem and the starting values,
  but is always independent of the time-step size $\tau$ and the time $t_n$.
\end{lemma}

\begin{proof}
The result follows from the proof of \cite[Lemma 2.3]{LiaoJiZhang:2020pfc}
with the help of \eqref{ieq: Linfty H2 embedding}.
\end{proof}

\section{Some discrete convolution inequalities}
\setcounter{equation}{0}

\subsection{Some properties of DOC-$\rmk$ kernels}

Our error analysis is closely related to the discrete convolution form \eqref{eq: DOC form PFC BDFk}, so we
need some detail properties of the DOC-$\rmk$ kernels $\theta^{(\rmk)}_{j}$ and
the associated discrete convolution inequalities.
At first, we have the following result.

\begin{lemma}\cite[Lemma 2.1]{LiaoTangZhou:2021bdf345}\label{lem: equvielant positive definite DOC}
The discrete kernels $b^{(\rmk)}_{j}$ in \eqref{def: BDF-k generating function} are positive (semi-)definite
if and only if the associated DOC-$\rmk$ kernels $\theta_{j}^{(\rmk)}$ in \eqref{def: DOC-Kernels}  are positive (semi-)definite.
\end{lemma}

Thanks to Lemma \ref{lem: BDF-k minimum eigenvalue} and Lemma \ref{lem: equvielant positive definite DOC}, the DOC-$\rmk$ kernels are positive definite. Moreover, we collect the decaying estimates in
\cite[Lemma 2.5]{LiaoTangZhou:2021bdf345} and obtain the following result.
\begin{lemma}\label{lem: two properties of DOC}
For $3\le \rmk\le 5$, the associated DOC-$\rmk$ kernels $\theta_{j}^{(\rmk)}$
defined in \eqref{def: DOC-Kernels} are positive definite
and satisfy the following decaying estimates
  $$\absb{\theta_{j}^{(\rmk)}}\le \frac{\rho_{\rmk}}{4}\braB{\frac{\rmk}{7}}^j\quad\text{for $j\ge0$,}$$
  where the constants $\rho_3=10/3$, $\rho_4=6$ and $\rho_5=96/5$.
\end{lemma}

To facilitate the convergence analysis, we present some discrete convolution
inequalities with respect to the DOC-$\rmk$ kernels $\theta_{j}^{(\rmk)}$.
For the BDF-$\rmk$ formula, consider the following matrices of order $m:=n-\rmk+1$
\begin{align}\label{matrix: Bk lower}
B_{\rmk,l}:=\left(
 \begin{array}{cccccc}
    b_0^{(\rmk)} & && &&\\
     b_1^{(\rmk)} &b_0^{(\rmk)} & &&&\\
      \vdots            &    \ddots          & \ddots & &    &\\
      b_{\rmk-1}^{(\rmk)}&\cdots &b_1^{(\rmk)}&b_0^{(\rmk)}&&\\
      &\ddots&\cdots&b_1^{(\rmk)}&b_0^{(\rmk)}&\\
      &&b_{\rmk-1}^{(\rmk)}&\cdots&b_1^{(\rmk)}&b_0^{(\rmk)}\\
  \end{array}
\right)_{m\times m}\quad\text{and}\quad B_{\rmk}:=B_{\rmk,l}+B_{\rmk,l}^T,
\end{align}
where $3\le \rmk \le 5$ and the index $n\ge\rmk$. Lemma \ref{lem: BDF-k minimum eigenvalue} says that
the real symmetric matrix $B_{\rmk}$ is positive definite. Moreover, introduce the matirces
\begin{align}\label{matrix: Theta-k lower}
\Theta_{\rmk,l}:=
\left(
\begin{array}{cccc}
\theta_{0}^{(\rmk)}  &                  &  & \\
\theta_{1}^{(\rmk)}  &\theta_{0}^{(\rmk)}  &  & \\
\vdots           &\vdots           &\ddots  &\\
\theta_{m-1}^{(\rmk)}&\theta_{m-2}^{(\rmk)}&\cdots  &\theta_{0}^{(\rmk)}  \\
\end{array}
\right)_{m\times m}\quad\text{and}\quad \Theta_{\rmk}:=\Theta_{\rmk,l}+\Theta_{\rmk,l}^T\,,
\end{align}
where the discrete kernels
$b_{j}^{(\rmk)}$ and $\theta_{j}^{(\rmk)}$ are defined by
\eqref{def: BDF-k generating function} and \eqref{def: DOC-Kernels}, respectively.
It follows from the discrete orthogonal
identity  \eqref{eq: orthogonal identity} that
\begin{align}\label{matrix: orthogonal identity}
\Theta_{\rmk,l}= B_{\rmk,l}^{-1},
\end{align}
and thus
\begin{align}\label{matrix: Theta-k formula via B-k}
\Theta_{\rmk}:=\Theta_{\rmk,l}+\Theta_{\rmk,l}^T
= B_{\rmk,l}^{-1}+( B_{\rmk,l}^{-1})^T
=(B_{\rmk,l}^{-1})^TB_{\rmk}B_{\rmk,l}^{-1}.
\end{align}
As stated in Lemma \ref{lem: two properties of DOC}, the real symmetric matrix $\Theta_{\rmk}$ is also positive definite.

\subsection{Eigenvalue estimates}

We present the following eigenvalue estimates of $B_{\rmk,l}^TB_{\rmk,l}$ and $\Theta_{\rmk}$
for any indexes $n\ge\rmk$.

\begin{lemma}\label{lem: tilde B2TB2-positiveDefinite}
There exists a positive constant  $\mathfrak{m}_{2\rmk}$
 such that $\lambda_{\max}\brab{ B_{\rmk,l}^T B_{\rmk,l}}\leq \mathfrak{m}_{2\rmk}$ for $3\le \rmk\le 5$.
\end{lemma}
\begin{proof}
For the matrix $B_{\rmk,l}$ in \eqref{matrix: Bk lower} of any order $m$,
 $B_{\rmk,l}^TB_{\rmk,l}$ is a real symmetric matrix no more than $(2\rmk-1)$ diagonals.
That is, each row of $B_{\rmk,l}^TB_{\rmk,l}$ has at most $(2\rmk-1)$ bounded elements computed from the BDF-$\rmk$ kernels $b_j^{(\rmk)}$.
\lan{The Gerschgorin's circle theorem implies that there is a finite bound  $\mathfrak{m}_{2\rmk}$
 such that $\lambda_{\max}\brab{ B_{\rmk,l}^T B_{\rmk,l}}\leq \mathfrak{m}_{2\rmk}$,
 also see a detailed computation of the constant $\mathfrak{m}_{22}$ in \cite[Lemma A2]{LiaoJiZhang:2020pfc} for the variable-step BDF2 fromula.}
\end{proof}

To avoid possible confusions,
we define the vector norm $\timenorm{\cdot}$ by
$\timenorm{\myvec{u}}:=\sqrt{\myvec{u}^T\myvec{u}}$ for any real vector $\myvec{u}$ and
the associated matrix norm
$\timenorm{A}:=\sqrt{\lambda_{\max}\brab{A^TA}}$.

\begin{lemma}\label{lem: Theta minimum eigenvalue}
The matrix $\Theta_{\rmk}$ in \eqref{matrix: Theta-k formula via B-k} satisfies
$\lambda_{\min}\brab{\Theta_{\rmk}}\geq \mathfrak{m}_{1\rmk}/\mathfrak{m}_{2\rmk}$ for $3\le \rmk\le 5$.
\end{lemma}
\begin{proof}
Lemma \ref{lem: BDF-k minimum eigenvalue} says that real symmetric  matrix $ B_{\rmk}$ is positive definite.
There exists a non-singular upper triangular matrix $U$ such that $ B_{\rmk}=U^TU$.
By using \eqref{matrix: Theta-k formula via B-k}, one gets
\begin{align*}
\myvec{v}^T\Theta_{\rmk}\myvec{v}=\myvec{v}^T(B_{\rmk,l}^{-1})^TB_{\rmk}B_{\rmk,l}^{-1}\myvec{v}
=(UB_{\rmk,l}^{-1}\myvec{v})^T UB_{\rmk,l}^{-1}\myvec{v}
=\timenorm{UB_{\rmk,l}^{-1}\myvec{v}}^2.
\end{align*}
Thus it follows that
\begin{align*}
\timenorm{\myvec{v}}^2=&\,\timenorm{B_{\rmk,l}U^{-1}UB_{\rmk,l}^{-1}\myvec{v}}^2
\leq\timenorm{B_{\rmk,l}U^{-1}}^2\timenorm{UB_{\rmk,l}^{-1}\myvec{v}}^2\\
\leq&\,\timenorm{ B_{\rmk,l}}^2\timenorm{U^{-1}}^2\myvec{v}^T\Theta_{\rmk}\myvec{v}
=\lambda_{\max}\brab{ B_{\rmk,l}^T B_{\rmk,l}}\lambda_{\max}\brab{B_{\rmk}^{-1}}\myvec{v}^T\Theta_{\rmk}\myvec{v}.
\end{align*}
Thus Lemmas \ref{lem: BDF-k minimum eigenvalue}
and \ref{lem: tilde B2TB2-positiveDefinite} yield the claimed result.
\end{proof}

\begin{lemma}\label{lem: Theta maximum eigenvalue}
There exists a positive constant  $\mathfrak{m}_{3\rmk}$
 such that $\lambda_{\max}\brab{\Theta_{\rmk}}\le \mathfrak{m}_{3\rmk}$ for $3\le \rmk\le 5$.
\end{lemma}
\begin{proof} The decaying properties of the DOC kernels $\theta_{j}^{(\rmk)}$ determine
the boundedness of the maximum eigenvalue of $(\Theta_{\rmk})_{m\times m}$.
For an arbitrary order $m$, Lemma \ref{lem: two properties of DOC} shows that
\begin{align*}
\mathfrak{R}_{m,j}:=&\,\sum_{\ell=1}^j\absb{\theta_{j-\ell}^{(\rmk)}}
+\sum_{\ell=j}^m\absb{\theta_{\ell-j}^{(\rmk)}}
\le\frac{\rho_{\rmk}}{4}\sum_{\ell=1}^j\braB{\frac{\rmk}{7}}^{j-\ell}
+\frac{\rho_{\rmk}}{4}\sum_{\ell=j}^m\braB{\frac{\rmk}{7}}^{\ell-j}
<\frac{7\rho_{\rmk}}{2(7-\rmk)}
\end{align*}
for $1\le j\le m$. One takes $\mathfrak{m}_{3\rmk}:=\frac{7\rho_{\rmk}}{2(7-\rmk)}$ such that
$\lambda_{\max}\brab{\Theta_{\rmk}}
\le\max_{1\le j\le m}\mathfrak{R}_{m,j}<\mathfrak{m}_{3\rmk}$
by the Gerschgorin's circle theorem. It completes the proof.
\end{proof}

\subsection{Discrete convolution inequalities}

The following lemmas describe some discrete convolution inequalities of Young-type.
Here and hereafter, we always denote
$\sum_{\ell,j}^{n,\ell}:=\sum_{\ell=\rmk}^n\sum_{j=\rmk}^\ell$
for the simplicity of presentation.

\begin{lemma}\label{lem:DOC quadr form Young inequ2}
For any $\varepsilon > 0$, any real sequence $\{v^{\ell}\}_{\ell=\rmk}^n$ and $\{w^{\ell}\}_{\ell=\rmk}^n$, it holds that
\begin{align*}
\sum_{\ell,j}^{n,\ell}\theta_{\ell-j}^{(\rmk)} v^jw^{\ell}
\le&\,\varepsilon\sum_{\ell,j}^{n,\ell}\theta_{\ell-j}^{(\rmk)} v^{\ell} v^j
+\frac{1}{2\mathfrak{m}_{1\rmk}\varepsilon}\sum_{\ell=\rmk}^{n}(w^{\ell})^2
\le\varepsilon\sum_{\ell,j}^{n,\ell}\theta_{\ell-j}^{(\rmk)} v^{\ell} v^j
+\frac{\mathfrak{m}_{2\rmk}}{\mathfrak{m}_{1\rmk}^2\varepsilon}
\sum_{\ell,j}^{n,\ell}\theta_{\ell-j}^{(\rmk)}w^jw^{\ell}.
\end{align*}
\end{lemma}
\begin{proof}
Let $\myvec{w}:=\brat{w^{\rmk},w^{\rmk+1},\cdots,w^n}^T$.
A similar proof of \cite[Lemma A.3]{LiaoJiZhang:2020pfc} gives
\begin{align}\label{lemproof:Young inequ2}
\sum_{\ell,j}^{n,\ell}\theta_{\ell-j}^{(\rmk)}v^j w^{\ell}
\le \varepsilon\sum_{\ell,j}^{n,\ell}\theta_{\ell-j}^{(\rmk)}v^j v^{\ell}
+\frac{1}{2\varepsilon}\myvec{w}^T B_{\rmk}^{-1}\myvec{w}
\quad \text{for any $\varepsilon > 0$}.
\end{align}
From the proof Lemma \ref{lem: Theta minimum eigenvalue}, we have
$B_{\rmk}^{-1}= U^{-1}\brat{U^{-1}}^T$ and then
\begin{align*}
\myvec{w}^T B_{\rmk}^{-1}\myvec{w}=&\,\myvec{w}^T U^{-1}
\brat{U^{-1}}^T\myvec{w}
=\timenorm{\brat{U^{-1}}^T \myvec{w}}^2\leq\timenorm{\brat{U^{-1}}^T}^2
\timenorm{\myvec{w}}^2\\
=&\,\lambda_{\max}\brab{(B_{\rmk})^{-1}}
\myvec{w}^T\myvec{w}\leq \frac{1}{\mathfrak{m}_{1\rmk}}\,\myvec{w}^T\myvec{w}
=\frac{1}{\mathfrak{m}_{1\rmk}}\sum_{\ell=\rmk}^{n}(w^{\ell})^2,
\end{align*}
where Lemma \ref{lem: BDF-k minimum eigenvalue} has been used.
Inserting it into \eqref{lemproof:Young inequ2}, we obtain the first claimed inequality.
The second inequality follows  immediately from Lemma  \ref{lem: Theta minimum eigenvalue},
which gives the minimum eigenvalue estimate of $\Theta_{\rmk}$. It completes the proof.
\end{proof}

\begin{lemma}\label{lem:DOC quadr form Young inequ-embedding}
For any $\varepsilon > 0$,  any real sequences $\{v^{\ell}\}_{\ell=\rmk}^n$ and $\{w^{\ell}\}_{\ell=\rmk}^n$, it holds that
\begin{align*}
\sum_{\ell,j}^{n,\ell}\theta_{\ell-j}^{(\rmk)} v^jw^{\ell}
\le&\,\varepsilon\sum_{\ell=\rmk}^n(v^{\ell})^2
+\frac{\mathfrak{m}_{3\rmk}}{4\mathfrak{m}_{1\rmk}\varepsilon}
\sum_{\ell=\rmk}^{n} (w^{\ell})^2
\le\varepsilon\sum_{\ell=\rmk}^n(v^{\ell})^2
+\frac{\mathfrak{m}_{2\rmk}\mathfrak{m}_{3\rmk}}{2\mathfrak{m}_{1\rmk}^2\varepsilon}
\sum_{\ell,j}^{n,\ell}\theta_{\ell-j}^{(\rmk)}w^jw^{\ell}.
\end{align*}
\end{lemma}

\begin{proof}Taking
$\varepsilon:=2\varepsilon_2/\mathfrak{m}_{3\rmk}$ in the first inequality of
Lemma \ref{lem:DOC quadr form Young inequ2} yields
\begin{align*}
\sum_{\ell,j}^{n,\ell}\theta_{\ell-j}^{(\rmk)} v^jw^{\ell}
\le&\,\frac{2\varepsilon_2}{\mathfrak{m}_{3\rmk}}\sum_{\ell,j}^{n,\ell}\theta_{\ell-j}^{(\rmk)} v^{\ell} v^j
+\frac{\mathfrak{m}_{3\rmk}}{4\mathfrak{m}_{1\rmk}\varepsilon_2}\sum_{\ell=\rmk}^{n}  (w^{\ell})^2\\
\le&\,\varepsilon_2\sum_{\ell=\rmk}^{n}(v^{\ell})^2
+\frac{\mathfrak{m}_{3\rmk}}{4\mathfrak{m}_{1\rmk}\varepsilon_2}
\sum_{\ell=\rmk}^{n}(w^{\ell})^2,
\end{align*}
where Lemma \ref{lem: Theta maximum eigenvalue} was used in the last inequality.
The first inequality is verified by choosing $\varepsilon_2:=\varepsilon$,
and the second one follows from Lemma  \ref{lem: Theta minimum eigenvalue}  immediately.
\end{proof}

\begin{lemma}\label{lem:inner quad ineq}
Let $v^n\in\mathbb{V}$ be a
	sequence of grid functions. For any constant $\varepsilon>0$,
\begin{align*}
\sum_{\ell,j}^{n,\ell}\theta_{\ell-j}^{(\rmk)}
\myinnerb{\Delta v^j,\Delta v^{\ell}}
\le&\,\varepsilon\sum_{\ell,j}^{n,\ell}\theta_{\ell-j}^{(\rmk)}
\myinnerb{\nabla\Delta v^j,\nabla\Delta v^{\ell}}
+\frac{8\mathfrak{m}_{2\rmk}^2}{\mathfrak{m}_{1\rmk}^5\varepsilon^2}\sum_{\ell=\rmk}^{n}\mynormb{v^\ell}^2.
\end{align*}
\end{lemma}

\begin{proof}
For any constant $\varepsilon_3>0$, we can apply
the second inequality of Lemma \ref{lem:DOC quadr form Young inequ2} with
$w^\ell:=-\nabla\Delta v^\ell$, $v^j:=\nabla v^j$ and $\varepsilon:=\mathfrak{m}_{2\rmk}
/\brat{\mathfrak{m}_{1\rmk}^2\varepsilon_3}$ and derive that
\begin{align*}
2\sum_{\ell,j}^{n,\ell}\theta_{\ell-j}^{(\rmk)}
\myinnerb{\Delta v^j,\Delta v^\ell}
\le&\,\frac{2\mathfrak{m}_{2\rmk}}
{\mathfrak{m}_{1\rmk}^2\varepsilon_3}
\sum_{\ell,j}^{n,\ell}\theta_{\ell-j}^{(\rmk)}
\myinnerb{\nabla v^j,\nabla v^\ell}
+2\varepsilon_3\sum_{\ell,j}^{n,\ell}\theta_{\ell-j}^{(\rmk)}
\myinnerb{\nabla\Delta v^j,\nabla\Delta v^\ell}.
\end{align*}
Similarly, by using the first inequality of 
Lemma \ref{lem:DOC quadr form Young inequ2}
with  $v^j:=-\Delta v^j$, $w^\ell:=v^\ell$  and the parameter $\varepsilon:=\varepsilon_3\mathfrak{m}_{1\rmk}^2
/(2\mathfrak{m}_{2\rmk})$, we can get
\begin{align*}
\frac{2\mathfrak{m}_{2\rmk}}{\mathfrak{m}_{1\rmk}^2\varepsilon_3}
\sum_{\ell,j}^{n,\ell}\theta_{\ell-j}^{(\rmk)}
\myinnerb{\nabla v^j,\nabla v^\ell}
\le&\,\sum_{\ell,j}^{n,\ell}\theta_{\ell-j}^{(\rmk)}
\myinnerb{\Delta v^j,\Delta v^\ell}
+\frac{2\mathfrak{m}_{2\rmk}^2}{\mathfrak{m}_{1\rmk}^5\varepsilon_3^2}
\sum_{\ell=\rmk}^{n}\mynormb{v^\ell}^2.
\end{align*}
We complete the proof by summing up the above two inequalities
and taking $\varepsilon_3:=\varepsilon/2$.
\end{proof}

\section{$L^2$ norm error estimate}
\setcounter{equation}{0}

Let $\xi^j:=D_{\rmk}\Phi(t_j)-\partial_t\Phi(t_j)$ be the local consistency error
of BDF-$\rmk$ formula at the time $t=t_j$.
Assume that the solution is regular in time for $t\ge t_{\rmk}$ such that
\begin{align}\label{ieq: time consisitency error}
\absb{\xi^j}\le C_\phi\tau^{\rmk}\max_{t_{\rmk}\le t\le T}\absb{\partial_t^{(\rmk+1)}\Phi(t)}
\le C_\phi\tau^{\rmk}\quad
\text{for $j\ge\rmk$}.
\end{align}
Then Lemma \ref{lem: two properties of DOC} yields
\begin{align}\label{ieq: BDF-k-global consistency}
\sum_{\ell=\rmk}^n\tau\mynormb{\Xi^\ell}\le C_\phi\tau^{\rmk+1} \sum_{\ell=\rmk}^n\sum_{j=\rmk}^{\ell}\absb{\theta_{\ell-j}^{(\rmk)}}
\le \frac{\rho_{\rmk}t_{n-\rmk+1}}{7-\rmk}C_\phi \tau^{\rmk}\quad\text{for $n\ge \rmk$,}
\end{align}
where the global time consistency error is defined by
\begin{align}\label{def: BDF-k-global consistency}
\Xi^\ell:=\sum_{j=\rmk}^\ell\theta_{\ell-j}^{(\rmk)}\xi^j\;\;\text{for $\ell\ge \rmk$.}
\end{align}

Note that,  the energy dissipation law \eqref{cont:energy dissipation}
of PFC model \eqref{cont: Problem-PFC} shows that $E[\Phi^n]\le E[\Phi(t_0)]$.
From the formulation \eqref{cont:free energy}, it is not difficult to see that
$\mynormb{\Phi^n}_{H^2}$ can be bounded by a time-independent constant. Applying the Sobolev embedding inequality, one has
\begin{align}\label{maximum bound exact solution}
\mynormb{\Phi^n}_{L^{\infty}}\le c_\Omega\mynormb{\Phi^n}_{H^2} 
\le c_2\quad\text{for $1\le n\le N$,}
\end{align}
where $c_2$ is dependent on the domain $\Omega$
and initial data $\Phi(t_0)$, but independent of the time $t_n$.

In the convergence analysis, set
$$c_3:=c_2^2+c_1c_2+c_1^2+\epsilon\quad
\text{and}\quad c_4:=250\mathfrak{m}_{2\rmk}^2/\mathfrak{m}_{1\rmk}^5
+2c_3^2\mathfrak{m}_{2\rmk}\mathfrak{m}_{3\rmk}/\mathfrak{m}_{1\rmk}^2,
$$
which may be dependent on the given data, the solution and
the starting values, but are always independent of 
the time-step size $\tau$ and the time $t_n$. Recall the following estimates on the starting values $\phi_{\mathrm{I}}^{(\rmk,j)}$ defined in
\eqref{Dis: initial effect -BDF formula}.
\begin{lemma}\cite[Lemma 2.6]{LiaoTangZhou:2021bdf345}\label{lem: Estimates of initial terms uI}
There exist some positive constants $c_{\mathrm{I},\rmk}>1$ such that
the starting values $\phi_{\mathrm{I}}^{(\rmk,j)}$ satisfy
\begin{align*}
\absb{\phi_{\mathrm{I}}^{(\rmk,j)}}\le \frac{c_{\mathrm{I},\rmk}\rho_{\rmk}}{8}
\braB{\frac{\rmk}{7}}^{j-\rmk}\sum_{\ell=1}^{\rmk-1}
\absb{\diff \phi^\ell}\quad\text{for $3\le \rmk\le 5$ and $j\ge\rmk$,}
\end{align*}
such that
\begin{align*}
\sum_{j=\rmk}^n\absb{\phi_{\mathrm{I}}^{(\rmk,j)}}\le
\frac{7c_{\mathrm{I},\rmk}\rho_{\rmk}}{8(7-\rmk)}
\sum_{\ell=1}^{\rmk-1}
\absb{\diff \phi^\ell}\quad\text{for $3\le \rmk\le 5$ and $n\ge\rmk$,}
\end{align*}
where the constants $\rho_{\rmk}$ are defined in Lemma \ref{lem: two properties of DOC}.
\end{lemma}

\begin{theorem}\label{thm:Convergence-Results}
Assume that the PFC problem \eqref{cont: Problem-PFC} has a solution $\Phi\in C^{\rmk+1}[0,T]$. 
If $\mathbf{A2}$ holds and the time-step size is small such that $\tau\le 1/(2c_4)$,
the numerical solution $\phi^n$ of the BDF-$\rmk$ implicit scheme \eqref{scheme: PFC BDFk}
is convergent in the $L^2$ norm,
\begin{align*}
\mynormb{\Phi^n-\phi^n}
\le \frac{7\rho_{\rmk}}{7-\rmk}\exp\brat{2c_4t_{n-\rmk+1}}&\,\Big(
c_{\mathrm{I},\rmk}\sum_{\ell=0}^{\rmk-1}\mynormb{\Phi^\ell-\phi^\ell}
+C_{\phi}t_{n-\rmk+1}\tau^{\rmk}\Big),
\quad\text{$\rmk\le n\le N$.}
\end{align*}
\end{theorem}

\begin{proof}
Let $e^n:=\Phi^n-\phi^n$ be the error
between the exact solution
and the numerical solution of the BDF-$\rmk$ implicit scheme \eqref{scheme: PFC BDFk}.
We have the following error equation
\begin{align}\label{Error-Equation}
D_{\rmk}e^n
&=\Delta\kbrab{(1+\Delta)^2e^n +f_{\phi}^ne^n}
+\xi^n\quad\text{for $\rmk\le n\le N$,}
\end{align}
where $\xi^n$ is defined by \eqref{ieq: time consisitency error} 
and $f_{\phi}^n :=\brat{\Phi^n}^2+\Phi^n\phi^n+\brat{\phi^n}^2-\epsilon$.
Thanks to the maximum norm solution estimates in Lemma \ref{lem:Bound-Solution} and \eqref{maximum bound exact solution},
one has
\begin{align}\label{ieq: max bound nonlinear}
\mynormb{f_{\phi}^n}_{\infty}\le c_2^2+c_1c_2+c_1^2+\epsilon=c_3.
\end{align}

Multiplying both sides of equation \eqref{Error-Equation} by $\tau\theta_{\ell-n}^{(\rmk)}$
and summing up $n$ from $n=\rmk$ to $\ell$, we apply the equality
\eqref{Dis: DOC action BDF formula Dk} with $v^j:=e^j$ to obtain
\begin{align}\label{Error-Equation-DOC}
\diff e^{\ell}
=-e_{\mathrm{I}}^{(\rmk,\ell)}+\tau\sum_{j=\rmk}^\ell
\theta_{\ell-j}^{(\rmk)}\Delta\kbrab{
(1+\Delta)^2e^j + f_{\phi}^j e^j}
+\tau\Xi^\ell\quad\text{for $\rmk\le \ell\le N$,}
\end{align}
where $\Xi^\ell$ is defined by 
\eqref{def: BDF-k-global consistency}
and $e_{\mathrm{I}}^{(\rmk,n)}$ represents the starting error effects
on the numerical solution at the time $t_n$
\begin{align}\label{def: starting error effects -BDF formula}
 e_{\mathrm{I}}^{(\rmk,n)}:=\sum_{\ell=1}^{\rmk-1}\diff e^\ell
\sum_{j=\rmk}^{n}\theta_{n-j}^{(\rmk)}b_{j-\ell}^{(\rmk)}\qquad\text{for $n\ge\rmk$.}
\end{align}
Making the inner product of \eqref{Error-Equation-DOC} with $2e^\ell$,
and summing up the superscript from $\rmk$ to $n$, we have the following equality
\begin{align}\label{Error-Equation-Inner}
\mynormb{e^n}^2-\mynormb{e^{\rmk-1}}^2
&\le-2\sum_{\ell=\rmk}^n\myinnerb{e_{\mathrm{I}}^{(\rmk,\ell)},e^{\ell}}+J^n +2\tau\sum_{\ell=\rmk}^n\myinnerb{\Xi^\ell,e^\ell}\quad\text{for $\rmk\le n\le N$,}
\end{align}
where the identity $2a(a-b)=a^2-b^2+(a-b)^2$ is used 
and $J^n$ is defined by
\begin{align}\label{Error-quadratic forms}
J^n:=&\,2\tau\sum_{\ell,j}^{n,\ell}\theta_{\ell-j}^{(\rmk)}
\myinnerb{e^j+2\Delta e^j+\Delta^2e^j+f_{\phi}^je^j,\Delta e^\ell}\nonumber\\
=&\,2\tau\sum_{\ell,j}^{n,\ell}\theta_{\ell-j}^{(\rmk)}
\kbra{\myinnerb{f_{\phi}^je^j+2\Delta e^j,\Delta e^\ell}-\myinnerb{\nabla e^j,\nabla e^\ell}
-\myinnerb{\nabla\Delta e^j,\nabla\Delta e^\ell}}.
\end{align}
Now we handle the quadratic form $J^n$.
By applying the second inequality of Lemma \ref{lem:DOC quadr form Young inequ-embedding}
with $v^j:=f_{\phi}^j e^j$, $w^\ell:=\Delta e^\ell$ and
$\varepsilon:=\mathfrak{m}_{2\rmk}\mathfrak{m}_{3\rmk}
/{\mathfrak{m}_{1\rmk}^2}$, one derives that
\begin{align*}
2\tau\sum_{\ell,j}^{n,\ell}\theta_{\ell-j}^{(\rmk)}
&\,\myinnerb{f_{\phi}^je^j+2\Delta e^j,\Delta e^\ell}
=2\tau\sum_{\ell,j}^{n,\ell}\theta_{\ell-j}^{(\rmk)}
\myinnerb{f_{\phi}^je^j,\Delta e^\ell}
+4\tau\sum_{\ell,j}^{n,\ell}\theta_{\ell-j}^{(\rmk)}
\myinnerb{\Delta e^j,\Delta e^\ell}\\
\le&\, \frac{2\mathfrak{m}_{2\rmk}\mathfrak{m}_{3\rmk}}
{\mathfrak{m}_{1\rmk}^2}
\sum_{\ell=\rmk}^{n}\tau\mynormb{f_{\phi}^\ell e^\ell}^2
+5\tau\sum_{\ell,j}^{n,\ell}\theta_{\ell-j}^{(\rmk)}
\myinnerb{\Delta e^j,\Delta e^\ell}\\
\le&\,\frac{2c_3^2\mathfrak{m}_{2\rmk}\mathfrak{m}_{3\rmk}}
{\mathfrak{m}_{1\rmk}^2}
\sum_{\ell=\rmk}^{n}\tau\mynormb{e^\ell}^2
+\frac{250\mathfrak{m}_{2\rmk}^2}{\mathfrak{m}_{1\rmk}^5}
\sum_{\ell=\rmk}^{n}\tau\mynormb{e^\ell}^2
+2\tau\sum_{\ell,j}^{n,\ell}\theta_{\ell-j}^{(\rmk)}
\myinnerb{\nabla\Delta e^j,\nabla\Delta e^\ell},
\end{align*}
where the maximum norm estimate \eqref{ieq: max bound nonlinear}, and
Lemma \ref{lem:inner quad ineq} with $v^j:=e^j$ and
$\varepsilon:=2/5$ were used in the second inequality.
Also, Lemma \ref{lem: two properties of DOC} implies that
$-\sum_{\ell,j}^{n,\ell}\theta_{\ell-j}^{(\rmk)}
\myinnerb{\nabla e^j,\nabla e^l}\le0$.
Then we obtain from \eqref{Error-quadratic forms} that
\begin{align*}
J^n\le c_4\sum_{\ell=\rmk}^{n}\tau\mynormb{e^\ell}^2.
\end{align*}
Therefore, it follows from \eqref{Error-Equation-Inner} that
\begin{align*}
\mynormb{e^n}^2
\le \mynormb{e^{\rmk-1}}^2+2\sum_{\ell=\rmk}^n\mynormb{e_{\mathrm{I}}^{(\rmk,\ell)}}\mynormb{e^{\ell}}
+c_4\sum_{\ell=\rmk}^n\tau   \mynormb{e^{\ell}}^2
+ 2\tau\sum_{\ell=\rmk}^n\mynormb{e^{\ell}}\mynormb{\Xi^\ell}
\quad\text{for $\rmk\le n\le N.$}
\end{align*}
Choosing some integer $n_0$ ($\rmk-1\le n_0 \le n$) such that
$\mynormb{e^{n_0}}=\max_{\rmk-1\le \ell \le n}\mynormb{e^{\ell}}$.
Taking $n:=n_0$ in the above inequality, one can obtain
\begin{align*}
\mynormb{e^{n_0}}\le \mynormb{e^{\rmk-1}}+2\sum_{\ell=\rmk}^{n_0}\mynormb{e_{\mathrm{I}}^{(\rmk,\ell)}}
+c_4\sum_{\ell=\rmk}^{n_0}\tau   \mynormb{e^{\ell}}+ 2\tau\sum_{\ell=\rmk}^{n_0}
\mynormb{\Xi^\ell}.
\end{align*}
By applying Lemma \ref{lem: Estimates of initial terms uI} to the starting term $e_{\mathrm{I}}^{(\rmk,\ell)}$
in \eqref{def: starting error effects -BDF formula}, one has
\begin{align*}
2\sum_{\ell=\rmk}^n\mynormb{e_{\mathrm{I}}^{(\rmk,\ell)}}
\leq \frac{7c_{\mathrm{I},\rmk}\rho_{\rmk}}{4(7-\rmk)}
\sum_{\ell=1}^{\rmk-1}\mynormb{\diff e^\ell}\quad\text{for $\rmk\le n\le N$}.
\end{align*}
Thus one gets
\begin{align*}
\mynormb{e^n}\le\mynormb{e^{n_0}}
&\le \frac{7c_{\mathrm{I},\rmk}\rho_{\rmk}}{2(7-\rmk)}
\sum_{\ell=0}^{\rmk-1}\mynormb{e^\ell}
+c_4\sum_{\ell=\rmk}^{n}\tau   \mynormb{e^{\ell}}
+ 2\tau\sum_{\ell=\rmk}^{n}\mynormb{\Xi^\ell}.
\end{align*}
Under the time-step constraint $\tau\le 1/(2c_4)$, we have
\begin{align*}
\mynormb{e^n}\le \frac{7c_{\mathrm{I},\rmk}\rho_{\rmk}}{7-\rmk}
\sum_{\ell=0}^{\rmk-1}\mynormb{e^\ell}
+2c_4\sum_{\ell=\rmk}^{n-1}\tau   \mynormb{e^{\ell}}
+ 4\tau\sum_{\ell=\rmk}^n\mynormb{\Xi^\ell}.
\end{align*}
By the standard discrete Gr\"onwall inequality together with the consistency 
estimate \eqref{ieq: BDF-k-global consistency}, 
one can obtain the claimed error estimate and complete the proof.
\end{proof}

\section{Numerical experiments}
\setcounter{equation}{0}

Some numerical experiments  are included to illustrate the efficiency of the BDF-$\rmk$ schemes
by the Fourier pseudo-spectral method in space. The resulting nonlinear algebraic systems are solved by fixed-point iterative methods with the termination error $10^{-12}$. Due to periodic boundary conditions, the fast Fourier transform can be applied for every iteration step. The sixth-order Gauss collocation method \cite{GongZhaoWang:2020HEQRK} is employed to initiate the numerical schemes such that the assumptions $\mathbf{A1}$ and $\mathbf{A2}$ would be reasonably fulfilled.

\begin{example}
We consider the exterior-forced PFC model $\partial_t \Phi=\Delta \mu +g(\mathbf{x},t)$ with the model parameter $\epsilon=0.02$,  which has an exact solution $\Phi=\cos(t)\sin(\frac{\pi}{2}x)\sin(\frac{\pi}{2}y)$.
\end{example}

\begin{table}[htb!]
	\begin{center}
		\caption{Numerical accuracy of BDF-$\rmk$ scheme 
			\eqref{scheme: PFC BDFk}.}
		\label{exam:PFC-RK-BDFk-Time-Error} \vspace*{0.3pt}
		\def\temptablewidth{0.92\textwidth}
		{\rule{\temptablewidth}{0.5pt}}
			\begin{tabular*}{\temptablewidth}
				{@{\extracolsep{\fill}}cccccccc}
				\multirow{2}{*}{$N$}  &\multirow{2}{*}{$\tau$}
				&\multicolumn{2}{c}{BDF$3$ scheme}  &\multicolumn{2}{c}{BDF$4$ scheme}  &\multicolumn{2}{c}{BDF$5$ scheme}\\
				\cline{3-4}          \cline{5-6}    \cline{7-8}    
			    &	&$e(N)$    &Order     &$e(N)$   &Order 
			    &$e(N)$   &Order  \\
				\midrule
				10   & 1.00e-01 &1.85e-04  &$-$   &1.19e-05   &$-$   
				&3.85e-06   &$-$ \\
				20   &5.00e-02  &2.42e-05  &2.94  &7.97e-07   &3.90  
				&9.53e-08   &5.34 \\
				40  &2.50e-02  &3.08e-06  &2.98  &5.14e-08   &3.95 
				&1.80e-09   &5.72 \\
				80  &1.25e-02  &3.71e-07  &3.05  &3.26e-09   &3.98 
				&3.86e-11   &5.54 \\
				160  &6.25e-03  &4.60e-08  &3.01  &2.05e-10   &4.00
				&1.16e-12   &5.06 \\
			\end{tabular*}
		{\rule{\temptablewidth}{0.5pt}}
	\end{center}
\end{table}

The domain $\Omega=(0,8)^2$ is divided into a $128 \times 128$ mesh
such that the temporal error dominates the spatial error in each run.
We solve the problem until time $T=1$.
The numerical result is tabulated in Table \ref{exam:PFC-RK-BDFk-Time-Error},
in which the discrete $L^2$ norm error $e(N):=\|\Phi(T)-\phi^N\|$ is recorded in each run
and the experimental  order is computed by
$\text{Order}\approx\log_2\bra{e(N)/e(2N)}$. It is observed that the BDF-$\rmk$
scheme is $\rmk$th-order accuracy in time.

\begin{figure}[htb!]
\begin{center}
\includegraphics[width=1.90in]{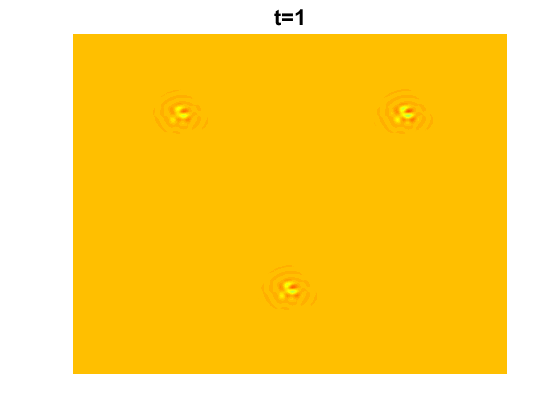}
\includegraphics[width=1.90in]{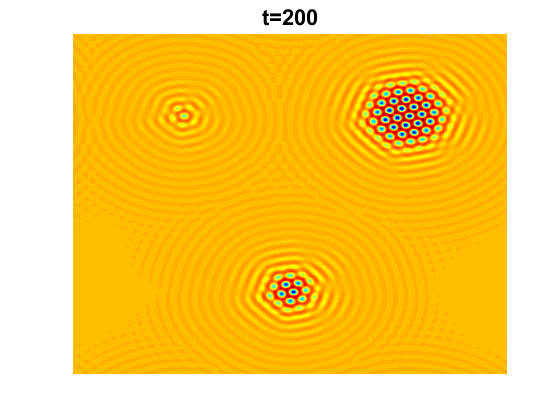}
\includegraphics[width=1.90in]{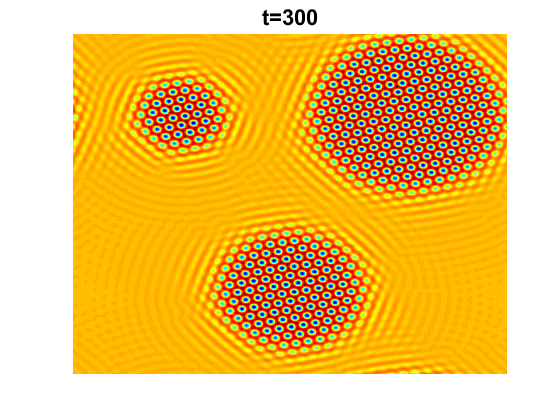}\\
\includegraphics[width=1.90in]{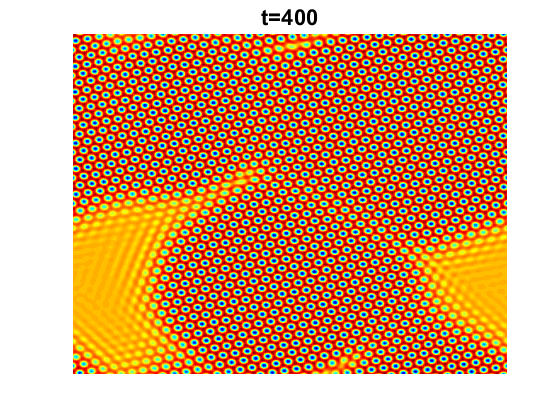}
\includegraphics[width=1.90in]{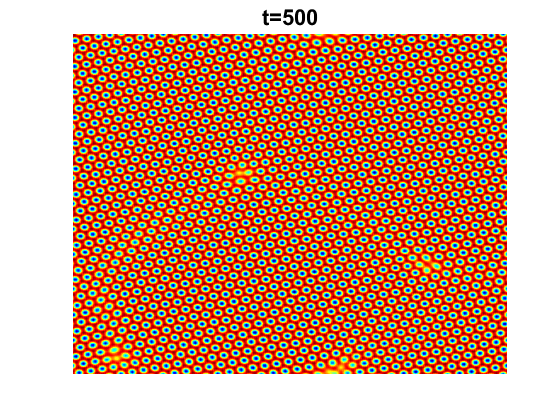}
\includegraphics[width=1.90in]{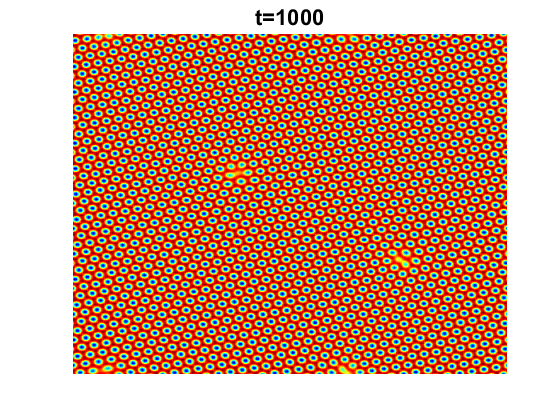}\\
\end{center}
\caption{The crystal growth process obtained at
$t=1$, $200$, $300$, $400$, $500$ and $1000$ by the BDF-5 scheme
(the BDF-3 and BDF-4 schemes generate similar profiles).}
\label{exam:BDF5 solution evolution}
\end{figure}

%

\begin{figure}[htb!]
\begin{center}
\includegraphics[width=2.6in]{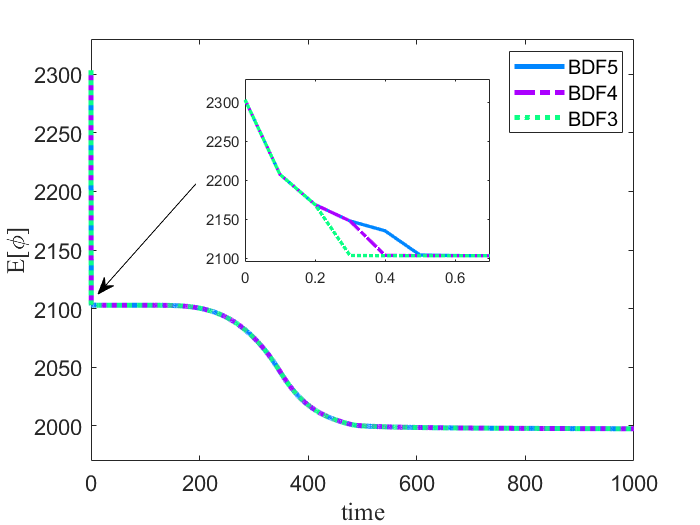}
\includegraphics[width=2.6in]{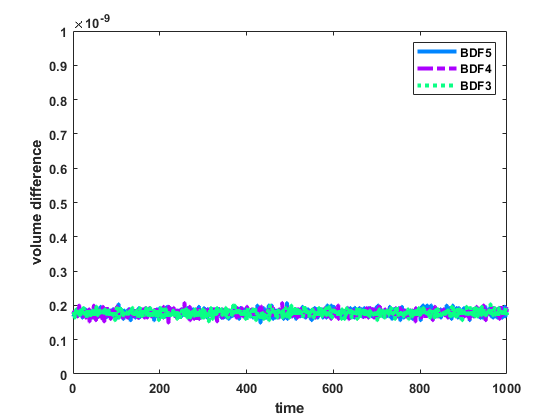}
\end{center}
\caption{Evolutions of original energy (left) and volume difference (right)}
\label{exam: Energy volume}
\end{figure}

\begin{example}
 We take the parameter $\epsilon=0.25$ and use a $256 \times 256$ uniform mesh to discretize the spatial domain $\Omega=(0,256)^{2}$. As seeds for nucleation, three random perturbations on the three small square patches are taken as $\Phi_{0}(\boldsymbol{x})=\bar{\Phi}+\operatorname{A \cdot  rand}(\boldsymbol{x})$, where the constant density $\bar{\Phi}=0.285, A$ is amplitude and the random numbers rand $(\cdot)$ are uniformly distributed in $(-1,1) .$ The centers of three pathes locate at  $(64,196)$, $(128,64)$ and $(196,196)$, with the corresponding amplitudes  $A=0.25$, $0.3$ and $0.35$, respectively. The length of each small square is set to 10. The solution is computed until the time $T=1000$ with a constant time step $\tau=0.1$.
\end{example}

The time evolutions of the phase variable are depicted in Figure \ref{exam:BDF5 solution evolution}. It is clear that the speed of moving interfaces is related to the initial amplitude A, the larger the amplitude A, the faster the crystal growth. And  three different crystal grains grow and become large enough to form grain boundaries eventually.
The discrete original  energy \eqref{def: discrete free energy} and the volume difference are shown in  Figure \ref{exam: Energy volume}. As predicted by our theory, the discrete volume is conservative (up to a tolerance  $10^{-9}$). It can be seen that the energy dissipates very fast at the early stage, and gradually slows down as the time escapes.


\appendix

\section{Proof of Lemma \ref{lem: Quadratic decomposition BDF-k}}
\label{section: quadratic decomposition BDF-k}
\setcounter{equation}{0}

This proof describes some quadratic decompositions of the following quantity
\begin{align*}
\mathfrak{B}_\rmk^n:=v_n\sum_{j=1}^nb_{n-j}^{(\rmk)}v_j
=b_{0}^{(\rmk)}v_n^2+b_{1}^{(\rmk)}v_nv_{n-1}+\cdots+b_{\rmk-1}^{(\rmk)}v_nv_{n-\rmk+1}\quad\text{for $n\ge \rmk$,}
\end{align*}
where the discrete BDF-$\rmk$ kernels $b_{j}^{(\rmk)}$ (vanish if $j\ge \rmk$)
are generated by \eqref{def: BDF-k generating function}, see Table \ref{table: BDF-k kernels}.
According to the derivations of BDF-$\rmk$ formulas, we use the difference operators $\delta_1^{m}v_n$ to find
\begin{align*}
\sum_{j=1}^nb_{n-j}^{(\rmk)}v_j=\sum_{m=1}^{\rmk}\frac{1}{m}\delta_1^{m-1}v_n
=v_n+\frac{1}{2}\delta_1v_n+\frac{1}{3}\delta_1^{2}v_n+\cdots+\frac{1}{\rmk}\delta_1^{\rmk-1}v_n.
\end{align*}
It implies that
\begin{align}\label{QuadDecompo: BDF-k quantity}
\mathfrak{B}_\rmk^n
=&\,\sum_{m=1}^{\rmk}\frac{1}{m}v_n\delta_1^{m-1}v_n
=\sum_{m=1}^{\rmk}\frac{1}{2m}J_{m-1}[v_n],
\quad\text{where}\quad J_m[v_n]:=2v_n\delta_1^mv_n.
\end{align}
Obviously, we have two trivial cases,
$$\mathfrak{B}_1^n=\frac{1}{2}J_{0}[v_n]=v_n^2$$
and
\begin{align*}
\mathfrak{B}_2^n=&\,\mathfrak{B}_1^n+\frac{1}{4}J_{1}[v_n]
=\frac{1}{4}\brab{v_n^2-v_{n-1}^2}+v_n^2+\frac{1}{4}(\delta_1v_n)^2.
\end{align*}
In general, we will handle $J_{\rmk}[v_n]$ and decompose $\mathfrak{B}_{\rmk+1}^n$ via the following equality,
\begin{align}\label{QuadDecompo: BDF-k recurive procedure}
\mathfrak{B}_{\rmk+1}^n=&\,\mathfrak{B}_{\rmk}^n+\frac{1}{2(\rmk+1)}J_{\rmk}[v_n]\quad\text{for $\rmk\ge2$.}
\end{align}
We prove Lemma \ref{lem: Quadratic decomposition BDF-k} for the cases $\rmk=3,4$ and $5$ in the subsequent subsections, respectively.
Our process for the quadratic decompositions includes the following three steps:
\begin{itemize}[itemindent=1cm]
	\item [\textbf{Step 1.}] Apply the identities $2a(a-b)=a^2-b^2+(a-b)^2$ and  $2b(a-b)=a^2-b^2-(a-b)^2$ 
	to decompose $J_m[v_n]=2v_n\delta_1^mv_n$ in \eqref{QuadDecompo: BDF-k quantity} into some quadratic terms,
	see \eqref{QuadDecompo: J2 procedure}, \eqref{QuadDecompo: J3 procedure} and \eqref{QuadDecompo: J4 procedure}.
	Then one can obtain the preliminary (not necessarily desired) quadratic decomposition of  
	$\mathfrak{B}_\rmk^n$ via the recurive formula in \eqref{QuadDecompo: BDF-k recurive procedure}.
	\item [\textbf{Step 2.}] Apply the inverse decomposition formulas, see  \eqref{QuadDecompo: inverse J2 tool},
	\eqref{QuadDecompo: inverse J3 tool} and \eqref{QuadDecompo: inverse J4 tool}, to 
	absorb some nonpositive quadratic terms into the nonnegative terms of high-order difference. 
	See \underline{the underlined parts} in this proof, we use the nonnegative terms $(\delta_1^2v_{n}+\alpha v_{n-1})^2$, $(\delta_1^3v_{n}+\beta \delta_1v_{n-1})^2$ and $(\delta_1^4v_{n}+\gamma \delta_1^2v_{n-1})^2$ to absorb the nonpositive terms $-(\delta_1v_{n-1})^2$,	$-(\delta_1^2v_{n-1})^2$ and  $-(\delta_1^3v_{n-1})^2$, respectively. Here $\alpha$, $\beta$ and $\gamma$ are constants. 
	\item [\textbf{Step 3.}] Repeat Step 2 untill the preliminary quadratic decomposition in Step 1 
	can be reformulated into a discrete gradient structure like \eqref{def: discrete gradient structure}.
\end{itemize}

\subsection{Decomposition for the BDF-3 formula}

Consider the case of $\rmk=2$. By noticing that
\begin{align}\label{QuadDecompo: J2 tool}
2v_{n-1}\delta_1^2v_{n}=v_n^2-2v_{n-1}^2+v_{n-2}^2-(\delta_1v_{n})^2-(\delta_1v_{n-1})^2,
\end{align}
one has
\begin{align}\label{QuadDecompo: J2 procedure}
J_{2}[v_n]
=2\delta_1v_n(\delta_1^2v_n)+2v_{n-1}\delta_1^2v_{n}
=v_n^2-2v_{n-1}^2+v_{n-2}^2-2(\delta_1v_{n-1})^2+(\delta_1^2v_n)^2.
\end{align}
Then  we obtain
\begin{align*}
\mathfrak{B}_{3}^n=&\,\mathfrak{B}_{2}^n+\frac{1}{6}J_{2}[v_n]\nonumber\\
=&\,\frac{17}{12}v_n^2-\frac{7}{12}v_{n-1}^2+\frac{1}{6}v_{n-2}^2
+\frac{1}{4}(\delta_1v_n)^2-\frac{1}{4}(\delta_1v_{n-1})^2
+\underline{\frac{1}{6}(\delta_1^2v_n)^2-\frac{1}{12}(\delta_1v_{n-1})^2}.
\end{align*}
We treat with the last two terms (the underlined part) as follows,
\begin{align*}
\widetilde{R}_{31}^n:=&\,\frac{1}{6}(\delta_1^2v_n)^2-\frac{1}{24}(\delta_1v_{n-1})^2\\
=&\,\frac{1}{6}(\delta_1^2v_n)^2+\frac{1}{12}v_{n-1}\delta_1^2v_{n}
-\frac{1}{24}v_n^2+\frac{1}{12}v_{n-1}^2-\frac{1}{24}v_{n-2}^2+\frac{1}{24}(\delta_1v_{n})^2\\
=&\,\frac{1}{6}(\delta_1^2v_n+\tfrac14v_{n-1})^2
-\frac{1}{24}v_n^2+\frac{7}{96}v_{n-1}^2-\frac{1}{24}v_{n-2}^2+\frac{1}{24}(\delta_1v_{n})^2,
\end{align*}
where the equality \eqref{QuadDecompo: J2 tool} was applied inversely, that is,
\begin{align}\label{QuadDecompo: inverse J2 tool}
-(\delta_1v_{n-1})^2=2v_{n-1}\delta_1^2v_{n}-v_n^2+2v_{n-1}^2-v_{n-2}^2+(\delta_1v_{n})^2.
\end{align}
Then we derive that
\begin{align*}
\mathfrak{B}_{3}^n
=&\,\kbra{\frac{37}{96}v_{n}^2-\frac{1}{8}v_{n-1}^2+\frac{7}{24}(\delta_1v_{n})^2}
-\kbra{\frac{37}{96}v_{n-1}^2-\frac{1}{8}v_{n-2}^2+\frac{7}{24}(\delta_1v_{n-1})^2}\nonumber\\
&\,+\frac{95}{96}v_n^2+\frac{1}{6}(\delta_1^2v_n+\tfrac14v_{n-1})^2.
\end{align*}
Let $\sigma_{L3}:=95/48$, and
 \begin{align}
\mathcal{G}_{3}[\vec{v}_n]:=&\,\frac{37}{96}v_{n}^2-\frac{1}{8}v_{n-1}^2+\frac{7}{24}(\delta_1v_{n})^2
=\frac{1}{6}v_{n}^2+\frac1{6}(\tfrac7{4}v_n-v_{n-1})^2\ge0,
\label{QuadDecompo: BDF-3 pesdoEnergy}\\
\mathcal{R}_{3}[\vec{v}_n]:=&\,\frac{1}{6}(\delta_1^2v_n+\tfrac14v_{n-1})^2.
\label{QuadDecompo: BDF-3 Remainder}
\end{align}
It follows that
\begin{align}\label{QuadDecompo: BDF-3 decomposition}
\mathfrak{B}_{3}^n
=&\,\mathcal{G}_{3}[\vec{v}_n]-\mathcal{G}_{3}[\vec{v}_{n-1}]+\frac{\sigma_{L3}}2v_n^2+\mathcal{R}_{3}[\vec{v}_n].
\end{align}
It confirms the claimed
decomposition \eqref{QuadDecompo: general quadratic decomposition}
for the case of $\rmk=3$.

\subsection{Decomposition for the BDF-4 formula}

Consider the case of $\rmk=3$. One follows the
derivations of \eqref{QuadDecompo: J2 procedure} to obtain
\begin{align}\label{QuadDecompo: J3 procedure}
J_{3}[v_n]=&\,2v_n\delta_1^3v_n
=J_{2}[\delta_1v_n]+2v_{n-1}\delta_1^2v_{n}-J_{2}[v_{n-1}]\nonumber\\
=&\,-3(\delta_1v_{n-1})^2+3(\delta_1v_{n-2})^2-3(\delta_1^2v_{n-1})^2+(\delta_1^3v_n)^2\nonumber\\
&\,+v_n^2-3v_{n-1}^2+3v_{n-2}^2-v_{n-3}^2,
\end{align}
where the equality \eqref{QuadDecompo: J2 tool} was also used.
Then using the quadratic decomposition \eqref{QuadDecompo: BDF-3 decomposition}
together with \eqref{QuadDecompo: BDF-3 pesdoEnergy} and
\eqref{QuadDecompo: BDF-3 Remainder} we obtain
\begin{align}\label{QuadDecompo: BDF-4 decomposition0}
\mathfrak{B}_{4}^n=&\,\mathfrak{B}_{3}^n+\frac{1}{8}J_{3}[v_n]
=\mathcal{G}_{3}[\vec{v}_n]-\mathcal{G}_{3}[\vec{v}_{n-1}]+R_4^n,
\end{align}
where, by combining similar terms,
\begin{align*}
R_4^n:=&\,\frac{95}{96}v_n^2+\mathcal{R}_{3}[\vec{v}_n]+\frac{1}{8}J_{3}[v_n]
=\frac{107}{96}v_n^2-\frac{3}{8}v_{n-1}^2+\frac{3}{8}v_{n-2}^2-\frac{1}{8}v_{n-3}^2
\\
&\,-\frac{3}{8}(\delta_1v_{n-1})^2+\frac{3}{8}(\delta_1v_{n-2})^2
 +\frac{1}{6}(\delta_1^2v_n+\tfrac14v_{n-1})^2
+\underline{\frac{1}{8}(\delta_1^3v_n)^2-\frac{3}{8}(\delta_1^2v_{n-1})^2}.
\end{align*}
We will handle the last two terms (the underlined part).
Noticing that
\begin{align}\label{QuadDecompo: J3 tool}
2\delta_1v_{n-1}(\delta_1^3v_{n})
=(\delta_1v_n)^2-2(\delta_1v_{n-1})^2+(\delta_1v_{n-2})^2-(\delta_1^2v_{n})^2-(\delta_1^2v_{n-1})^2,
\end{align}
or, inversely,
\begin{align}\label{QuadDecompo: inverse J3 tool}
-(\delta_1^2v_{n-1})^2=2\delta_1v_{n-1}\delta_1^3v_{n}
-(\delta_1v_n)^2+2(\delta_1v_{n-1})^2-(\delta_1v_{n-2})^2+(\delta_1^2v_{n})^2,
\end{align}
one can derive that
\begin{align*}
\widetilde{R}_{41}^n:=&\,\frac{1}{8}(\delta_1^3v_n)^2-\frac{3}{16}(\delta_1^2v_{n-1})^2\\
=&\frac{1}{8}(\delta_1^3v_n+\tfrac3{2}\delta_1v_{n-1})^2
-\frac{3}{16}(\delta_1v_n)^2
+\frac{3}{32}(\delta_1v_{n-1})^2-\frac{3}{16}(\delta_1v_{n-2})^2
+\frac{3}{16}(\delta_1^2v_{n})^2.
\end{align*}
Inserting it into the above expression of $R_4^n$, one has
\begin{align*}
R_4^n
=&\,\frac{107}{96}v_n^2-\frac{3}{8}v_{n-1}^2+\frac{3}{8}v_{n-2}^2-\frac{1}{8}v_{n-3}^2\\
&\, +\frac{3}{16}(\delta_1^2v_{n})^2
-\frac{3}{16}(\delta_1^2v_{n-1})^2+\frac{1}{8}(\delta_1^3v_n+\tfrac3{2}\delta_1v_{n-1})^2\\
&\,-\frac{3}{16}(\delta_1v_n)^2-\frac{9}{64}(\delta_1v_{n-1})^2+\frac{3}{16}(\delta_1v_{n-2})^2
+\underline{\frac{1}{6}(\delta_1^2v_n+\tfrac14v_{n-1})^2-\frac{9}{64}(\delta_1v_{n-1})^2}.
\end{align*}
Now we handle the last two terms (the underlined part) by applying \eqref{QuadDecompo: inverse J2 tool} as follows,
\begin{align*}
\widetilde{R}_{42}^n:=&\,\frac{1}{6}(\delta_1^2v_n+\tfrac14v_{n-1})^2-\frac{9}{64}(\delta_1v_{n-1})^2\\
=&\,\frac{1}{6}(\delta_1^2v_n+\tfrac{35}{32}v_{n-1})^2
-\frac{9}{64}v_n^2+\frac{189}{32\cdot64}v_{n-1}^2-\frac{9}{64}v_{n-2}^2+\frac{9}{64}(\delta_1v_{n})^2.
\end{align*}
Inserting it into the above expression of $R_4^n$, one gets
\begin{align*}
R_4^n
=&\,\kbra{\frac{355}{2048}v_{n}^2-\frac{7}{64}v_{n-1}^2+\frac{1}{8}v_{n-2}^2}
-\kbra{\frac{355}{2048}v_{n-1}^2-\frac{7}{64}v_{n-2}^2+\frac{1}{8}v_{n-3}^2}\\
&\,-\frac{3}{64}(\delta_1v_n)^2+\frac{3}{64}(\delta_1v_{n-1})^2
-\frac{3}{16}(\delta_1v_{n-1})^2+\frac{3}{16}(\delta_1v_{n-2})^2
+\frac{3}{16}(\delta_1^2v_{n})^2\\
&\,-\frac{3}{16}(\delta_1^2v_{n-1})^2+\frac{4919}{6144}v_n^2+\frac{1}{8}(\delta_1^3v_n+\tfrac3{2}\delta_1v_{n-1})^2
+\frac{1}{6}(\delta_1^2v_n+\tfrac{35}{32}v_{n-1})^2.
\end{align*}
Inserting it into the equality \eqref{QuadDecompo: BDF-4 decomposition0}, one gets the desired decomposition
\begin{align}\label{QuadDecompo: BDF-4 decomposition}
\mathfrak{B}_{4}^n=\mathcal{G}_{4}[\vec{v}_n]-\mathcal{G}_{4}[\vec{v}_{n-1}]
+\frac{\sigma_{L4}}2v_n^2+\mathcal{R}_{4}[\vec{v}_n],
\end{align}
where $\sigma_{L4}=\frac{4919}{3072}$, the functionals $\mathcal{G}_{4}$ and $\mathcal{R}_{4}$ are defined by
     \begin{align}
\mathcal{G}_{4}[\vec{v}_n]:=&\,\frac{3433}{6144}v_{n}^2-\frac{15}{64}v_{n-1}^2+\frac{1}{8}v_{n-2}^2
+\frac{47}{192}(\delta_1v_{n})^2-\frac{3}{16}(\delta_1v_{n-1})^2+\frac{3}{16}(\delta_1^2v_{n})^2,
\label{QuadDecompo: BDF-4 pesdoEnergy}\\
\mathcal{R}_{4}[\vec{v}_n]:=&\,\frac{1}{8}(\delta_1^3v_n+\tfrac3{2}\delta_1v_{n-1})^2
+\frac{1}{6}(\delta_1^2v_n+\tfrac{35}{32}v_{n-1})^2.
\label{QuadDecompo: BDF-4 Remainder}
\end{align}
It confirms the claimed decomposition \eqref{QuadDecompo: general quadratic decomposition} for the case of $\rmk=4$,
because the quadratic functional $\mathcal{G}_{4}$ is non-negative, that is,
\begin{align*}
\mathcal{G}_{4}[\vec{v}_n]
=&\,\frac{3433}{6144}v_{n}^2-\frac{3}{16}v_n^2-\frac{3}{64}v_{n-1}^2+\frac{1}{8}v_{n-2}^2
+\frac{119}{192}(\delta_1v_{n})^2+\frac{3}{8}v_{n-2}(\delta_1v_{n})\\
=&\,\frac{13627}{43008}v_{n}^2+\frac{7}{24}(\tfrac{65}{56}v_n-v_{n-1})^2
+\frac{1}{8}(\tfrac3{2}\delta_1v_{n}+v_{n-2})^2\ge0.
\end{align*}

\subsection{Decomposition for the BDF-5 formula}

Consider the case of $\rmk=4$. By using \eqref{QuadDecompo: J2 tool} and \eqref{QuadDecompo: J2 procedure}, one has
\begin{align*}
2v_{n-1}\delta_1^3v_{n}
=&\,v_n^2-3v_{n-1}^2+3v_{n-2}^2-v_{n-3}^2
-(\delta_1v_{n})^2-(\delta_1v_{n-1})^2
+2(\delta_1v_{n-2})^2-(\delta_1^2v_{n-1})^2.
\end{align*}
Then we can follow the derivations of \eqref{QuadDecompo: J3 procedure} to obtain
\begin{align}\label{QuadDecompo: J4 procedure}
J_{4}[v_n]=&\,2v_n\delta_1^4v_n
=J_{3}[\delta_1v_n]+2v_{n-1}\delta_1^3v_{n}-J_{3}[v_{n-1}]\nonumber\\
=&\,v_n^2-4v_{n-1}^2+6v_{n-2}^2-4v_{n-3}^2+v_{n-4}^2
-4(\delta_1v_{n-1})^2+8(\delta_1v_{n-2})^2-4(\delta_1v_{n-3})^2\nonumber\\
&\,-4(\delta_1^2v_{n-1})^2+6(\delta_1^2v_{n-2})^2-4(\delta_1^3v_{n-1})^2+(\delta_1^4v_n)^2.
\end{align}
Then using \eqref{QuadDecompo: BDF-4 decomposition} together with \eqref{QuadDecompo: BDF-4 pesdoEnergy} and
\eqref{QuadDecompo: BDF-4 Remainder}, we obtain
\begin{align}\label{QuadDecompo: BDF-5 decomposition0}
\mathfrak{B}_{5}^n=&\,\mathfrak{B}_{4}^n+\frac{1}{10}J_{4}[v_n]
=\mathcal{G}_{4}[\vec{v}_n]-\mathcal{G}_{4}[\vec{v}_{n-1}]+R_{5}^n,
\end{align}
where, by combining similar terms,
\begin{align}\label{QuadDecompo: BDF-5 Remaider0}
R_{5}^n=&\,\frac{4919}{6144}v_n^2+\mathcal{R}_{4}[\vec{v}_n]+\frac{1}{10}J_{4}[v_n]
=\frac{4919}{6144}v_n^2+\frac{1}{8}(\delta_1^3v_n+\tfrac3{2}\delta_1v_{n-1})^2
+\frac{1}{6}(\delta_1^2v_n+\tfrac{35}{32}v_{n-1})^2
\nonumber\\
&\,+\frac{1}{10}v_n^2-\frac{2}{5}v_{n-1}^2+\frac{3}{5}v_{n-2}^2-\frac{2}{5}v_{n-3}^2+\frac{1}{10}v_{n-4}^2
-\frac{2}{5}(\delta_1v_{n-1})^2+\frac{4}{5}(\delta_1v_{n-2})^2\nonumber\\
&\,-\frac{2}{5}(\delta_1v_{n-3})^2-\frac{2}{5}(\delta_1^2v_{n-1})^2+\frac{3}{5}(\delta_1^2v_{n-2})^2
\underline{-\frac{2}{5}(\delta_1^3v_{n-1})^2+\frac{1}{10}(\delta_1^4v_n)^2}.
\end{align}
Noticing that
\begin{align}\label{QuadDecompo: J4 tool}
2\delta_1^2v_{n-1}(\delta_1^4v_{n})
=(\delta_1^2v_n)^2-2(\delta_1^2v_{n-1})^2+(\delta_1^2v_{n-2})^2-(\delta_1^3v_{n})^2-(\delta_1^3v_{n-1})^2,
\end{align}
or, inversely,
\begin{align}\label{QuadDecompo: inverse J4 tool}
-(\delta_1^3v_{n-1})^2=2\delta_1^2v_{n-1}(\delta_1^4v_{n})
-(\delta_1^2v_n)^2+2(\delta_1^2v_{n-1})^2-(\delta_1^2v_{n-2})^2+(\delta_1^3v_{n})^2,
\end{align}
we handle the last two terms (the underlined part) in \eqref{QuadDecompo: BDF-5 Remaider0} as follows
\begin{align*}
\widetilde{R}_{51}^n:=&\,\frac{1}{10}(\delta_1^4v_n)^2-\frac{1}{5}(\delta_1^3v_{n-1})^2-\frac{1}{5}(\delta_1^3v_{n-1})^2\\
=&\,\frac{1}{10}(\delta_1^4v_n+2\delta_1^2v_{n-1})^2
-\frac{1}{5}(\delta_1^2v_n)^2-\frac{1}{5}(\delta_1^2v_{n-2})^2+\frac{1}{5}(\delta_1^3v_{n})^2-\frac{1}{5}(\delta_1^3v_{n-1})^2.
\end{align*}
Then it follows from \eqref{QuadDecompo: BDF-5 Remaider0} that
\begin{align}\label{QuadDecompo: BDF-5 Remaider1}
R_{5}^n=&\,
\frac{4919}{6144}v_n^2+\frac{1}{10}v_n^2-\frac{2}{5}v_{n-1}^2+\frac{3}{5}v_{n-2}^2-\frac{2}{5}v_{n-3}^2
+\frac{1}{10}v_{n-4}^2\nonumber\\
&\,-\frac{2}{5}(\delta_1v_{n-1})^2+\frac{4}{5}(\delta_1v_{n-2})^2
-\frac{2}{5}(\delta_1v_{n-3})^2+\frac{2}{5}(\delta_1^2v_{n-2})^2
\nonumber\\
&\,+\frac{1}{5}(\delta_1^3v_{n})^2-\frac{1}{5}(\delta_1^3v_{n-1})^2
+\frac{1}{10}(\delta_1^4v_n+2\delta_1^2v_{n-1})^2
+\frac{1}{6}(\delta_1^2v_n+\tfrac{35}{32}v_{n-1})^2
\nonumber\\
&\,-\frac{1}{5}(\delta_1^2v_n)^2-\frac{3}{10}(\delta_1^2v_{n-1})^2
+\underline{\frac{1}{8}(\delta_1^3v_n+\tfrac3{2}\delta_1v_{n-1})^2-\frac{1}{10}(\delta_1^2v_{n-1})^2}.
\end{align}
By using \eqref{QuadDecompo: inverse J3 tool}, we treat with the last two terms (the underlined part) by
\begin{align*}
\widetilde{R}_{52}^n:=&\,\frac{1}{8}(\delta_1^3v_n+\tfrac3{2}\delta_1v_{n-1})^2-\frac{1}{10}(\delta_1^2v_{n-1})^2\\
=&\,\frac{1}{8}(\delta_1^3v_n+\tfrac{23}{10}\delta_1v_{n-1})^2
-\frac{1}{10}(\delta_1v_n)^2-\frac{9}{50}(\delta_1v_{n-1})^2
-\frac{1}{10}(\delta_1v_{n-2})^2+\frac{1}{10}(\delta_1^2v_{n})^2.
\end{align*}
Inserting it into \eqref{QuadDecompo: BDF-5 Remaider1}, we obtain
\begin{align}\label{QuadDecompo: BDF-5 Remaider2}
R_{5}^n=&\,
\frac{4919}{6144}v_n^2+\frac{1}{10}v_n^2-\frac{2}{5}v_{n-1}^2+\frac{3}{5}v_{n-2}^2-\frac{2}{5}v_{n-3}^2
+\frac{1}{10}v_{n-4}^2-\frac{1}{10}(\delta_1v_n)^2\nonumber\\
&\,-\frac{39}{100}(\delta_1v_{n-1})^2
+\frac{3}{10}(\delta_1v_{n-2})^2+\frac{2}{5}(\delta_1v_{n-2})^2-\frac{2}{5}(\delta_1v_{n-3})^2
-\frac{1}{10}(\delta_1^2v_n)^2\nonumber\\
&\,-\frac{3}{10}(\delta_1^2v_{n-1})^2+\frac{2}{5}(\delta_1^2v_{n-2})^2
+\frac{1}{5}(\delta_1^3v_{n})^2-\frac{1}{5}(\delta_1^3v_{n-1})^2
+\frac{1}{10}(\delta_1^4v_n+2\delta_1^2v_{n-1})^2\nonumber\\
&\,+\frac{1}{8}(\delta_1^3v_n+\tfrac{23}{10}\delta_1v_{n-1})^2
+\underline{\frac{1}{6}(\delta_1^2v_n+\tfrac{35}{32}v_{n-1})^2-\frac{19}{100}(\delta_1v_{n-1})^2}.
\end{align}
Furthermore, one can apply \eqref{QuadDecompo: inverse J2 tool} to get
\begin{align*}
\widetilde{R}_{53}^n:=&\,\frac{1}{6}(\delta_1^2v_n+\tfrac{35}{32}v_{n-1})^2-\frac{19}{100}(\delta_1v_{n-1})^2\\
=&\,\frac{1}{6}(\delta_1^2v_n+\tfrac{1787}{800}v_{n-1})^2
-\frac{19}{100}v_n^2-\frac{10089}{40000}v_{n-1}^2-\frac{19}{100}v_{n-2}^2+\frac{19}{100}(\delta_1v_{n})^2.
\end{align*}
Then we can derive from \eqref{QuadDecompo: BDF-5 Remaider2} that
\begin{align}\label{QuadDecompo: BDF-5 Remaider3}
R_{5}^n
=&\,\widetilde{\mathcal{G}}_{5}[\vec{v}_n]-\widetilde{\mathcal{G}}_{5}[\vec{v}_{n-1}]
+\frac{\sigma_{L5}}{2}v_n^2+\mathcal{R}_{5}[\vec{v}_n],
\end{align}
where the constant $\sigma_{L5}:=\frac{646631}{1920000}$, the functionals $\widetilde{\mathcal{G}}_{5}$ and $\mathcal{R}_{5}^n$ are defined by
\begin{align*}
\widetilde{\mathcal{G}}_{5}[\vec{v}_n]:=&\,\frac{21689}{40000}v_n^2-\frac{11}{100}v_{n-1}^2+\frac{3}{10}v_{n-2}^2
-\frac{1}{10}v_{n-3}^2+\frac{9}{100}(\delta_1v_n)^2-\frac{3}{10}(\delta_1v_{n-1})^2\nonumber\\
&\,+\frac{2}{5}(\delta_1v_{n-2})^2-\frac{1}{10}(\delta_1^2v_n)^2-\frac{2}{5}(\delta_1^2v_{n-1})^2
+\frac{1}{5}(\delta_1^3v_{n})^2,\nonumber\\
\mathcal{R}_{5}[\vec{v}_n]:=&\,
\frac{1}{10}(\delta_1^4v_n+2\delta_1^2v_{n-1})^2
+\frac{1}{8}(\delta_1^3v_n+\tfrac{23}{10}\delta_1v_{n-1})^2
+\frac{1}{6}(\delta_1^2v_n+\tfrac{1787}{800}v_{n-1})^2.
\end{align*}
Return to \eqref{QuadDecompo: BDF-5 decomposition0} and one gets
claimed decomposition \eqref{QuadDecompo: general quadratic decomposition} for $\rmk=5$,
\begin{align}\label{QuadDecompo: BDF-5 decomposition}
\mathfrak{B}_{5}^n=&\,\mathcal{G}_{5}[\vec{v}_n]-\mathcal{G}_{5}[\vec{v}_{n-1}]
+\frac{\sigma_{L5}}{2}v_n^2+\mathcal{R}_{5}[\vec{v}_n],
\end{align}
where, by using \eqref{QuadDecompo: BDF-4 pesdoEnergy},
\begin{align}
\mathcal{G}_{5}[\vec{v}_n]:=&\,\mathcal{G}_{4}[\vec{v}_n]+\widetilde{\mathcal{G}}_{5}[\vec{v}_n]\nonumber\\
=&\,\frac{4227769}{3840000}v_{n}^2-\frac{551}{1600}v_{n-1}^2+\frac{17}{40}v_{n-2}^2
-\frac{1}{10}v_{n-3}^2+\frac{1607}{4800}(\delta_1v_{n})^2\nonumber\\
&\,-\frac{39}{80}(\delta_1v_{n-1})^2+\frac{2}{5}(\delta_1v_{n-2})^2
+\frac{7}{80}(\delta_1^2v_{n})^2-\frac{2}{5}(\delta_1^2v_{n-1})^2
+\frac{1}{5}(\delta_1^3v_{n})^2.\label{QuadDecompo: BDF-5 pesdoEnergy}
\end{align}
By following the treatment of $\mathcal{G}_{4}$ in the above subsection, it is not difficult to find that
\begin{align*}
\mathcal{G}_{5}[\vec{v}_n]
=&\,\frac{1198850903}{1678080000}v_{n}^2
+\frac{437}{900}(\tfrac{4931}{6992}v_{n}-v_{n-1})^2\nonumber\\
&\,+\frac{9}{40}(\tfrac{23}{18}\delta_1v_{n}+v_{n-2})^2
+\frac{1}{10}(2\delta_1v_{n}+2v_{n-2}-v_{n-3})^2\ge0.
\end{align*}
The proof of Lemma \ref{lem: Quadratic decomposition BDF-k} is completed.


\end{document}